\documentclass[11pt, a4paper]{article}
\usepackage{graphicx}
\usepackage{amssymb,stmaryrd,amsmath}

 \input colordvi

     \def\White#1{\Color{0 0 0 0}{#1}}

\newtheorem{theorem}{Theorem}

\newtheorem{proposition}[theorem]{Proposition}
\newtheorem{lemma}[theorem]{Lemma}

\newcommand{\proof}{\noindent{\bf Proof \ }}
\def\qed{\hfill \mbox{$\Box$} \vspace{3mm}}

  \def\Req{\mbox{\ $=_{\scriptscriptstyle{\cal R}}$\ }}
  \def\Seq{\mbox{\ $=_{\scriptscriptstyle{\cal S}}$\ }}
  \def\R{\mbox{${\cal R}$}}
  \def\Q{\mbox{${\cal Q}$}}

  \def\A{\mbox{${\cal A}$}}
  \def\B{\mbox{${\cal B}$}}
  \def\Bprev{\mbox{$B_{\rev}^+$}}
  \def\Rer{\mbox{${\cal R}_{e,r}$}}
  \def\Aer{\mbox{${\cal A}_{e,r}$}}
  \def\Rdr{\mbox{${\cal R}_{d,r}$}}
  
  \def\C{\mbox{${\cal C}$}}
  
  \def\revstorelR{\,\mbox{$\curvearrowright_{\scriptscriptstyle{\cal R}}$}}
  
  \def\onerevstorelR{\,\mbox{$\curvearrowright\!\!\!\!\!\!{\scriptscriptstyle{_1}}_{\scriptscriptstyle{\ \cal R}}$}}
  \def\onerevstorelchiR{\,\mbox{$\curvearrowright\!\!\!\!\!\!{\scriptscriptstyle{_1}}_{\scriptscriptstyle{\ \chi({\cal R})}}$}}
  
  \def\revsto{\,\mbox{$\curvearrowright$ }}
  \def\revstoer{\,\mbox{$\curvearrowright_{e,r}$ }}
  \def\onerevsto{\,\mbox{$\curvearrowright\!\!\!\!\!\!{\scriptscriptstyle{_1}}\ \ $}}
  \def\gap{\mbox{\texttt{GAP}}}
  \def\mod{\mbox{\texttt{mod}}\,}
  \def\rev{\mbox{\texttt{rev}}}

  \def\phimt{\mbox{$\stackrel{\varphi}{\mapsto}$}}
  \def\bjmt{\mbox{$\stackrel{\alpha_{\mathbf j}}{\mapsto}$}}

  \def\bt{\mbox{$\alpha_{\mathbf t}$}}

\def\bop{\mbox{$\binom{p}{0}$}}
\def\boq{\mbox{$\binom{q}{0}$}}

\def\bpq{\mbox{$\binom{q}{p+1}^{-1} \binom{q}{p}$}}

\def\bqp{\mbox{$\bop^{-1} \boq^{-1} 1^{-1} 2^{-1} 321 \boq \bop$}}
\def\boi{\mbox{$\langle 21 \rangle^{i} ( \langle 21 \rangle^{i-1} )^{-1}$}}
\def\bzii{\mbox{$\left( b_0^{(i)} \right)^{-1}$}}
\def\bzi{\mbox{ $b_0^{(i)}$ }}

\title{Garside structure for the braid group of $G(e,e,r)$}
\author{David Bessis \& Ruth Corran}
\date{}

\begin{document}
\maketitle

\begin{abstract}
We give a new presentation of the braid group $B$ of the complex
reflection group $G(e,e,r)$ which is positive and homogeneous,
and for which the generators
map to reflections in the corresponding complex reflection group.
We show that this presentation gives rise to a Garside structure 
for $B$ with Garside element a kind of generalised Coxeter element, 
and hence obtain solutions to the word and conjugacy problems for $B$. 
\end{abstract}


\section*{Introduction}

The work in this paper arose from some questions posed in~\cite{bmr}.
In that paper, the authors associate a braid group to each complex reflection group.
For all but six exceptional irreducible complex reflection groups, they exhibit particular presentations for these braid groups.
These presentations have certain remarkable properties, including that the relations are positive and
homogeneous, and that the generators are ``generators of the monodromy'' (natural analogues of
reflection in braid groups).
However, many questions about the braid groups remain open. For example, solutions to the word and conjugacy problems
are not known in all cases.

We consider the family of imprimitive reflection groups $G(e,e,r)$, where $e$ and $r$ are positive integral parameters.
This family contains two infinite families of real reflection groups: $G(2,2,r)$ is the reflection group
of type $D_r$, whereas $G(e,e,2)$ is the dihedral group $I_2(e)$. We denote by $B(e,e,r)$ the associated
braid group.

Except when $e=2$ or $r=2$, Brou\'e-Malle-Rouquier's presentation for $B(e,e,r)$ does not
have good algorithmic properties: 
the main obstruction is that the monoid defined  by their
positive homogeneous presentation does not embed in the group. This negative answer
to Question 2.28 in \cite{bmr} was 
first observed in the second author's PhD thesis
(\cite{corranphd}) -- in the last section of the present
article, we prove that an even stronger obstruction holds.

This obstruction led us to search for a new presentation for $B(e,e,r)$.

Our new presentation has more generators than the original presentation; our generators are conjugates
of the original generators, thus are ``generators of the monodromy''.
It is positive and homogeneous.
 For $e=2$ or $r=2$,
it coincides with the ``dual'' presentation introduced in \cite{dualmonoid}. Our main result, Theorem \ref{main}, 
states that the monoid defined by our presentation is a Garside monoid.

Garside monoids were introduced in different ways 
in~\cite{bdm} and~\cite{dehpa}. We use the approach of~\cite{dehpa}
in this article. A 
monoid $M$ defined by a homogeneous presentation is \emph{Garside} 
if it is cancellative, the partial orders on $M$ defined by 
left and right division are lattices, and $M$ possesses a special 
\emph{Garside element}, for which the sets of 
left divisors and right divisors are equal and generate $M$.
The group of fractions of a Garside monoid is a \emph{Garside group},
which is the group defined by the same presentation 
(considered now as a group presentation) as the original monoid.
Garside monoids embed in their groups of fractions;
Garside groups have explicitly solvable word~\cite{dehpa} and conjugacy~\cite{picantin:conjugacy} 
problems.
Thus our Theorem \ref{main} provides the first known solutions to the word and conjugacy
problems in $B(e,e,r)$.

This article is in four sections. 
We begin by exhibiting a 
positive homogeneous presentation, defined in terms of two parameters
$e$ and $r$. 
In the second section, we show that the monoid defined by this presentation has
a \emph{complete presentation} (see~\cite{dehornoy:complete}),
and thus that it is a Garside monoid. Hence the group defined by 
the presentation is a Garside group. The third section is devoted to 
showing that this group is indeed the braid group of the complex reflection 
group $G(e,e,r)$. In the last section, we show that, with the original generating set of \cite{bmr}, one cannot
write a finite positive homogeneous presentation satisfying the embedding property.

{\flushleft \bf Note.} The proofs presented in this article rely on complex computations, including some case-by-case verifications
performed using the computer algebra system \gap. After completing this first version, we realised that computer calculations could be
avoided with a more geometric approach. While we will present this new approach in a forthcoming sequel,
we decided to post this article with our first proof, since the generality of the techniques used may be applicable to
other contexts. We also postpone explanation of the origin of our new presentation,
since this is not required for the algebraic approach used here, and is
better understood in the geometric setting of the sequel.

\goodbreak


\section{The presentation}

\subsection{Statement of the presentation}
\label{statement}

Here we give a presentation which we will show in Section 2 to be 
a presentation of the braid group of the reflection group $G(e,e,r)$.
For notational reasons, we will use $r=n+1$ throughout.
The generators are $\A = \A_0 \cup \A_1$ where 
$$\A_0 = \{a_{pq} \mid 0 \leq p,q  < n\} \quad \mbox{ and }
\A_1 = \{a_{p}^{(i)} \mid 0 \leq p < n, \ 0 \leq i < e \}.$$ 
{\quotation{\small \noindent We show in Theorem~\ref{main} that these generators map to reflections in the 
 reflection group $G(e,e,r)$, which can be thought of as the group of monomial
 $ r \times r$ matrices whose non-zero entries firstly must come from  
 $\{\exp\left(\frac{2 \pi i}{e}\right) \mid 0 \leq i < e \}$, and secondly, must multiply 
 to give 1. To give a little meaning to  the generators \A\ above, and why they are 
 so denoted, we describe their images in $G(e,e,r)$.
 We will label the entries in matrices over the indices $\{0,1,\ldots, n\}$.
 Let $E_{ij}$ denote 
 the $r \times r$ matrix with 0's everywhere but 1 in the $(i,j)$ place.
 Define a matrix $M_{ij}(a,b) := I_r - E_{ii} - E_{jj} + a E_{ij} + b E_{ji}$
 (which is like a permutation matrix corresponding to the transposition
 $(i \ j)$, but with the off-diagonal entries replaced by $a$ and $b$).
 Let $\zeta = \exp\left(\frac{2 \pi}{e}\right)$;
 then the matrices corresponding to our generators are:
 for $0 \leq p < q <n$,
 $$
 \begin{array}{c}
 a_{pq} \stackrel{\nu}{\longmapsto} M_{pq}(1,1) \quad \mbox{ and } \quad 
 a_{qp} \stackrel{\nu}{\longmapsto} 
 M_{pq}(\zeta, \zeta^{-1});
 \end{array}$$
 and for $0 \leq p < n$ and $0 \leq i < e$, 
 $$a_{p}^{(i)} \stackrel{\nu}{\longmapsto} 
 M_{pn}(\zeta^i, \zeta^{-i}).$$
}}

There are six types of relations, they are given modulo $e$ where 
appropriate. 
The sequence $(a_1,\ldots,a_k)$ of distinct integers \mod $n$ 
is said to be $n$-\emph{cyclically ordered} 
if the sequence of roots of unity $(\mu^{a_1}, \ldots, \mu^{a_k})$ 
has strictly increasing arguments, where $\mu$ is the $n$th 
root of unity $\exp \left( \frac{2 \pi}{n} \right)$. 
For example, if $0 \leq p < q < r < n$ then 
certainly $(p,q,r)$ is  $n$-cyclically ordered; 
so are $(q,r,p)$ and $(r,p,q)$; and these three are 
the only $n$-cyclically ordered sequences on these numbers.

\noindent
For every $n$-cyclically ordered sequence $(p,q,r,s)$ there are relations:
$$\begin{array}{ lcl}
 (\R_1) &  a_{pq} a_{rs} = a_{rs} a_{pq}    &  \mbox{ and } \\
 (\R_2) &  a_{ps} a_{qr} = a_{qr} a_{ps} .   &     
\end{array}$$
For every $n$-cyclically ordered sequence $(p,q,r)$ and for every $i$ ($0 
\leq i < e$) there are relations:
$$\begin{array}{lcl}
  (\R_3) &   a_{pq} a_{qr}  =  a_{qr} a_{pr} =  a_{pr} a_{pq}  &   \\
  (\R_4) &     a_{pq} a_r^{(i)} =    a_r^{(i)}   a_{pq} 
.                  &
\end{array}$$
For every pair $0 \leq p , q < n$ and for every $i$ ($0 \leq i < e$) 
there are relations:
$$\begin{array}{lcl}
  (\R_5) &   a_{pq} a_{q}^{(j)}  =  a_{q}^{(j)} a_{p}^{(i)} =  
a_{p}^{(i)} a_{pq}  & \\
\end{array}$$
where 
$$j = \left\{ \begin{array}{ll}
     i & \mbox{if } p<q,  \mbox{ and}\\
     i+1 \ \ \mod e \  & \mbox{otherwise}.
\end{array} \right.$$ 
Finally for every $p$ ($0 \leq p < n$) there are relations:
$$ (\R_6) \quad  \quad 
a_p^{(0)} a_p^{(e-1)} 
= \cdots = a_p^{(i)} a_{p}^{(i-1)}  = \cdots =  a_p^{(2)} a_p^{(1)}  
= a_p^{(1)} a_p^{(0)}.$$  
We will write $B^+$ for the monoid defined by this presentation.

\subsubsection*{$B^+$ is isomorphic to its reverse}

Let $\rev$ be the involution on $\A^*$
defined by $\rev(a_1 \cdots a_f) = a_f \cdots a_1$. 
Recall that the monoid $B^+$ is defined by the presentation
$\langle \A \mid \R \rangle$ where \A\ and \R\ are the 
generators and relations respectively defined at the beginning of
Section~\ref{statement}.
Define \Bprev\ to be the monoid defined by the presentation
$\langle \A \mid \rev(\R) \rangle$. We call this monoid the reverse of $B^+$.
Clearly $\rev$ defines an anti-isomorphism 
$B^+ {\longrightarrow} \Bprev$.
We will show that they are in fact isomorphic.

Let $\chi$ be any map: $\A \longrightarrow \B$. Then
$\chi$ extends to a homomorphism between the free monoids
$\A^*$ and $\B^*$. 
For any set of relations \R\ over \A, 
write $\chi(\R)$ for the set 
$\{\chi(u) = \chi(v) \mid (u = v) \in \R \}$
of relations over \B. 
Then $\chi$ defines a homomorphism between the monoids 
$M_1 = \langle \A \mid \R \rangle$ and $M_2 = \langle \B \mid \chi(\R) \rangle$.
In particular, $\chi(M_1) = M_2$ precisely when $\chi(\A) = \B$, 
and $\chi(M_1) \cong M_1$ precisely when $\chi:\A \rightarrow \B$ is one-one.
Thus $\chi:M_1 \rightarrow M_2$ is an isomorphism whenever $\chi: \A \rightarrow \B$
is a bijection.

Define an operator $\chi$ on our set of generators \A\ by 
$$
\label{chidef}
\chi(a_{pq}) := a_{q'p'} \quad \quad \mbox{ and } \quad \quad 
\chi(a_{p}^{(i)}) := \left\{
\begin{array}{ll}
a_{p'}^{(i')} & \mbox{ if } p \neq 0, \mbox{ and }\\
a_{0}^{(i'-1)} & \mbox{ otherwise,}
\end{array} \right.
$$ where  $p' := n-p$, $q' = n-q$ and $i' = e-i$. 
Clearly $\chi$ is a permutation of \A,  
so by the previous paragraph, we have an isomorphism
\begin{equation}
\label{chisom}
\chi : \B^+ \stackrel{\sim}{\longrightarrow} \langle \A \mid \chi(\R) \rangle.
\end{equation}
Moreover,

\begin{lemma}
\label{anti}
\emph{$\chi(\R) = \rev(\R)$}
\end{lemma}

\proof
Since $\chi$ and \rev\ are bijections on words,
it suffices to show that $\rev(\chi(\R)) \subseteq \R$,
as since \rev\ is an involution this implies $\chi(\R) \subseteq \rev(\R)$,
which in turn by bijectivity of $\chi$, and finiteness of \R, 
ensures that $\chi(\R) = \rev(\R)$.

We will use $p'$ and $i'$ to denote $n-p$ and $e-i$ respectively,
with  context (subscript \emph{vs} exponent) making it clear which is intended.
Observe that any sequence $(p_1,\ldots,p_f)$ is cyclically ordered 
precisely when $(p'_f,\ldots,p'_1)$  is cyclically ordered, which is
precisely when $(p'_1,\ldots,p'_f)$ is reverse cyclically ordered. 
In the rest of this proof, $x<y<z$ will mean cyclically ordered, 
and $0<x<y<n$ will mean linearly (usually) ordered.

Consider the relations of type $\R_1$ and $\R_2$; for $p<q<r<s$ we have
$$\begin{array}{cl}
\rev(\chi(a_{pq} a_{rs} = a_{rs} a_{pq}) ) 
= \rev(a_{q'p'} a_{s'r'} = a_{s'r'} a_{q'p'})
= (a_{s'r'} a_{q'p'} = a_{q'p'} a_{s'r'}) & \mbox{ and} \\
\rev(\chi(a_{ps} a_{qr} = a_{qr} a_{ps}) ) 
= \rev(a_{s'p'} a_{r'q'} = a_{r'q'} a_{s'p'})
= (a_{r'q'} a_{s'p'} = a_{s'p'} a_{r'q'}) 
\end{array}
$$
which are relations of type $\R_1$ and $\R_2$ respectively also, 
with $s'<r'<q'<p'$.

Next consider the typical relation of type $\R_3$ for $p<q<r$:
$$\begin{array}{rcl}
\rev(\chi(a_{pq} a_{qr} = a_{qr} a_{pr} = a_{pr} a_{pq}) ) 
&= &\rev(a_{q'p'} a_{r'q'} = a_{r'q'} a_{r'p'} = a_{r'p'} a_{q'p'}) \\
&= &\ \ \ \ \ (a_{q'p'} a_{r'p'} = a_{r'p'} a_{r'q'} = a_{r'q'} a_{q'p'}),
\end{array}$$
which is a relation of type $\R_3$ with $r'<q'<p'$. Looking at a relation of type $\R_4$:
$$\rev(\chi(a_{pq} a_r^{(i)} = a_r^{(i)} a_{pq}) ) 
= \rev(a_{q'p'} a_{r'}^{(j)} =  a_{r'}^{(j)} a_{q'p'}) 
= (a_{q'p'} a_{r'}^{(j)} = a_{r'}^{(j)} a_{q'p'}),$$
where $j = i'-1$ if $r = 0$ and $j=i'$ otherwise;  this is a relation
of type $\R_4$ also, with $r'<q'<p'$.

We will consider the image of relations of type $\R_5$ in four cases. Suppose
$0<p<q<n$. Then we have $0<q'<p'<n$, and:
$$\begin{array}{rcl}
\rev(\chi(a_{0q} a_{q}^{(i)} = a_{q}^{(i)} a_{0}^{(i)} = a_{0}^{(i)} a_{0q}) ) 
&= &\rev(a_{q'0} a_{q'}^{(i')} = a_{q'}^{(i')} a_{0}^{(i'+1)} = a_{0}^{(i'+1)} a_{q'0}) \\
&= &\ \ \ \ \ (a_{q'0} a_{0}^{(i'+1)} = a_{0}^{(i'+1)} a_{q'}^{(i')} = a_{q'}^{(i')} a_{q'0}),\\
\ \\
\rev(\chi(a_{pq} a_{q}^{(i)} = a_{q}^{(i)} a_{p}^{(i)} = a_{p}^{(i)} a_{pq}) ) 
&= &\rev(a_{q'p'} a_{q'}^{(i')} = a_{q'}^{(i')} a_{p'}^{(i')} = a_{p'}^{(i')} a_{q'p'}) \\
&= &\ \ \ \ \ (a_{q'p'} a_{p'}^{(i')} = a_{p'}^{(i')} a_{q'}^{(i')} = a_{q'}^{(i')} a_{q'p'}),\\
\ \\
\rev(\chi(a_{q0} a_{0}^{(i+1)} = a_{0}^{(i+1)} a_{q}^{(i)} = a_{q}^{(i)} a_{q0}) ) 
&= &\rev(a_{0q'} a_{0}^{(i')} = a_{0}^{(i')} a_{q'}^{(i')} = a_{q'}^{(i')} a_{0q'}) \\
&= &\ \ \ \ \ (a_{0q'} a_{q'}^{(i')} = a_{q'}^{(i')} a_{0}^{(i')} = a_{0}^{(i')} a_{0q'}),\\
\ \\
\rev(\chi(a_{qp} a_{p}^{(i+1)} = a_{p}^{(i+1)} a_{q}^{(i)} = a_{q}^{(i)} a_{qp}) ) 
&= &\rev(a_{p'q'} a_{p'}^{(i'-1)} = a_{p'}^{(i'-1)} a_{q'}^{(i')} = a_{q'}^{(i')} a_{p'q'}) \\
&= &\ \ \ \ \ (a_{p'q'} a_{q'}^{(i')} = a_{q'}^{(i')} a_{p}^{(i'-1)} = a_{p'}^{(i'-1)} a_{p'q'}),
\end{array}$$
where in the third and fourth calculations we use the fact that 
$(i+1)' = e-i-1 = i'-1$, and so $\chi(a_0^{(i+1)}) = a_0^{(i')}$. Each of these
calculations results in a relation of type $\R_5$ also, although permuting the
types.

Finally we consider the relation of type $\R_6$: for every $0 < p < n$ we have
$$\begin{array}{rl}
\rev(\chi(a_{p}^{(0)} a_{p}^{(e-1)} = &\cdots =  a_{p}^{(i)} a_{p}^{(i-1)} =
 \cdots = a_{p}^{(2)} a_{p}^{(1)} = a_{p}^{(1)} a_{p}^{(0)} )) \\
 = &\rev(a_{p'}^{(0)} a_{p'}^{(1)} = \cdots =  a_{p'}^{(i')} a_{p'}^{(i'+1)} =
\cdots = a_{p'}^{(e-2)} a_{p'}^{(e-1)} = a_{p'}^{(e-1)} a_{p'}^{(0)} )\\
 = &(a_{p'}^{(0)} a_{p'}^{(e-1)} = a_{p'}^{(e-1)} a_{p'}^{(e-2)} = \cdots =  
a_{p'}^{(i'+1)}a_{p'}^{(i')}  =
\cdots = a_{p'}^{(1)} a_{p'}^{(0)} ),
\end{array}$$
which is clearly a relation of type $\R_6$ also. The case where the subscript is $0$
is similar, with the exponents shifted by 1; also resulting in a relation of type $\R_6$.
Thus we have that $\rev(\chi(\R)) \subseteq \R$, and hence $\chi(\R) = \rev(\R)$, as desired.
\qed

Substituting the equality thus obtained into the isomorphism~(\ref{chisom}) above, we obtain
\begin{theorem}
\label{revisom}
The monoid $B^+$ is isomorphic to its reverse.
\end{theorem}

\subsubsection*{Some consequences of the relations}

In this subsection we observe a number of easily derived equations which are 
consequences of the relations.
For every $n$-cyclically ordered sequence $(p,q,r,s)$ we have
$$
\begin{array}{lc }
  ( \C_1)  & a_{pr} a_{pq} a_{rs}  = a_{qs} a_{ps} a_{qr}    \\
  ( \C_2)  &  a_{ps} a_{pq} a_{qr}  a_s^{(1)} a_s^{(0)}  = a_{rq}   
a_{rp} a_{rs} a_q^{(1)}  a_q^{(0)}    
\end{array}
$$
For every $n$-cyclically ordered sequence $(p,q,r)$ and for every $i$ ($0 
\leq i < e$) we have
$$
\begin{array}{lc }
  ( \C_3)  & \quad  a_{pr} a_{pq} a_{r}^{(1)}  a_{r}^{(0)}  = a_{qp} 
a_{qr}  a_{p}^{(1)}  a_{p}^{(0)}   \quad    \\
  ( \C_4)  &  \quad   a_{pr} a_{pq} a_r^{(k)}  = a_q^{(i)}   a_{qr} 
a_p^{(j)}      \quad  
\end{array}
$$
where
$$ j = \left\{
\begin{array}{ll }
     i & \mbox{ if } p < q, \mbox{ and }   \\
     i-1 \ \mod e & \mbox{ otherwise;}  
\end{array} \right.
\quad \quad k = \left\{
\begin{array}{ ll}
     i &   \mbox{ if } q < r, \mbox{ and }   \\
     i+1 \ \mod e&   \mbox{ otherwise.} 
\end{array} \right.$$
And lastly, for each pair $p, q$  and for each pair $i, j$ we have
$$\begin{array}{lc}
( \C_5) & \quad \quad  \quad  \  \ 
a_{pq}  a_q^{(1)} a_q^{(0)} = a_{qp}  a_p^{(1)} a_p^{(0)} \    \quad 
\quad   \\
( \C_6) & \quad \quad  \quad  \  \ 
a_p^{(i)} a_p^{(i-1)} a_{pq} = a_q^{(j)} a_q^{(j-1)} a_{qp}  \    
\quad \quad   \\
\end{array}$$
In fact, if we denote by $t_p$ the 
word $a_p^{(1)} a_p^{(0)} $, we have that
$a_{pq} t_q = t_q a_{qp} = a_{qp}  t_p = t_p a_{pq} $,
from which $( \C_5)$  and $( \C_6)$  follow directly.

Call the presentation consisting of the original generating set,
along with relations $(\R_1), \ldots, (\R_6)$ together
with their consequences $(\C_1), \ldots, (C_6)$ the 
\emph{augmented presentation}. Since the added relations are consequences
of the original relations, the monoid and group presented by it are the
same as for the original presentation.
We observe that for every pair of generators
$x$ and $y$ there exists exactly one relation of the form $xu=yv$ in the
augmented presentation.


\section{Garside Structure}

We want to show that $B^+$ is a Garside monoid with Garside element
$\Delta$, the least common multiple of the generators. 
To do so, we will show that it is cancellative, that 
the partial orders on $M$ defined by 
left and right division are lattices, and that the sets of 
left divisors and right divisors of $\Delta$ are equal and generate $M$.


\subsection{Completeness}


A positive presentation $\cal P$ is said to be \emph{right complemented} if 
firstly, for every generator $x$, $\cal P$ contains no relations
of the form $xu = xv$ for any words $u$ and $v$, and 
secondly, for all generators $x \neq y$, there is exactly one relation
of the form $xu = yv$. Thus, by inspection, the augmented presentation for $B^+$ 
defined above can be seen to be right complemented. 

Completeness is introduced in \cite{dehornoy:complete}.
If a monoid possesses a presentation
which is both right complete and right complemented
then in particular the monoid is left cancellative (6.2 of~\cite{dehornoy:complete})
and the poset defined by left division is a lattice
(6.10 of~\cite{dehornoy:complete}).
We will show in Theorem~\ref{complete} that it the augmented presentation for
$B^+$ right complete.
Invoking Theorem~\ref{revisom} we will then have that $B^+$ is both 
left and right cancellative, and that the posets defined by left and
right division are lattices, as desired.

\subsubsection*{Word reversing}

Completeness is defined in~\cite{dehornoy:complete} via a technique
called \emph{word reversing}. We introduce this in the context of
a right complemented (positive) presentation $\langle \A \mid \R \rangle$. 
We begin by defining word reversing on generators: for
any generators $x$ and $y$, write 
$$x^{-1} y \onerevstorelR 
\left\{
\begin{array}{ll}
\varepsilon & \mbox{if $x=y$, and}\\
u v^{-1}&\mbox{if $xu=yv$ is a relation from $\R$.}
\end{array}
\right.$$
We extend in the obvious way to arbitrary words over $X = \A \cup \A^{-1}$:
for any $w, w'$ from $X^*$, write
$w x^{-1} y w' \onerevstorelR w z w'$ where $x^{-1} y \onerevstorelR z$.
We say that \emph{$w$ reverses to $w'$ in one step}.
Finally, for any $w, w'$ from $X^*$, 
we write $w \revstorelR w'$ if there exists a sequence of words 
$w \equiv w_0, w_1, \ldots, w_f \equiv w'$ from $X^*$ such that
$w_{i-1} \onerevstorelR w_i$ for each $i \in \{1,\ldots f\}$. In this case
we say that \emph{$w$ reverses to $w'$} (modulo \R). (The notation $w \equiv w'$
indicates that the words $w$ and $w'$ are letterwise identical.)

While reversing a word of the form $s^{-1} t$ where $s$ and $t$ are generators 
gives a unique result
(due to the presentation being complemented), in general a word may be
reversed in a single step to any of a number of other words: for example
$s^{-1} s t^{-1} t$ reverses in one step to $t^{-1} t$ and to $s^{-1} s$.
However there is confluence (again due to being complemented): if 
$w \revstorelR w'$ and $w \revstorelR w''$ 
where $w'$ and $w''$ are both fully word reversed 
(that is, no further reversings are possible),
then $w' \equiv w''$. 
Thus the particular order chosen to fully reverse a word is not important:
the end result is unique.

Clearly, a word which is fully word reversed is one in which
there is no occurrence of a subword of the form $s^{-1} t$ for 
generators $s$ and $t$. Thus these are precisely the words of the form
$uv^{-1}$ where $u$ and $v$ are words over \A. (For $v \equiv a_{i_1} \cdots a_{i_f}$ 
a word over \A,  we write $v^{-1}$ to denote the word $a_{i_f}^{-1} \cdots a_{i_f}^{-1}$.)

\subsubsection*{Definition of completeness}

We say that a triple of generators $(x,y,z)$ over \A\ satisfies the
\emph{completion condition modulo \R\ }if 
\begin{center}
\fbox{\parbox{135mm}{
\vspace{0mm}{
whenever $x^{-1} z z^{-1} y \revstorelR u v^{-1}$ 
for $u$ and $v$ positive words, 
we have  $(xu)^{-1} yv \revstorelR \varepsilon$}
\vspace{1mm}}}
\end{center}
where $\varepsilon$ represents the empty word. It is shown in~\cite{dehornoy:complete} 
that a presentation $\langle \A \mid \R \rangle$ is right complete if the 
completion condition holds modulo \R\ for every triple of generators over \A.
Until the end of this section we use the notation $\A_{e,r}$ and $\R_{e,r}$
to refer to the generators and augmented relations respectively of the 
braid group of type $G(e,e,r)$. (Thus $r=n+1$ in the presentation given
in Section~\ref{statement}.) For ease on the eye, we write
$\revstoer$ for $\curvearrowright_{\scriptscriptstyle{\cal R}_{e,r}}$.

We programmed \texttt{GAP}~\cite{gap} to take as input a set (or subset) 
of generators $G$ and a set of relations $R$, to calculate for each triple
$(x,y,z)$ over $G$ whether the above condition holds. The procedure outputs
to a file the calculations at each one-reversing step. The code to do this,
as well as the output calculations we will refer to below can be found at
\verb"http//www.math.jussieu.fr/~corran".

For example, this can be used to check, for fixed values of $e$ and $r$, 
that the presentation $\langle A_{e,r} \mid \R_{e,r} \rangle$ is complete.
In order to prove the presentation is complete for general $e$ and $r$
we will show that checking the completion condition modulo $\R_{e,r}$ 
for a triple $(x,y,z)$ over $\A_{e,r}$ can be reduced to 
checking the completion condition for another triple, 
this time modulo small values of $e$ and $r$, 
which then is done with the \gap\ code 
(see Propositions~\ref{gencalcs} and~\ref{allr}).

In Table~\ref{onesteps} (found on page~\pageref{onesteps}) 
we list all the one-step reversings of words of the form $x^{-1} y$
where $x$ and $y$ are generators; that is, 
for a pair of generators $x$ and $y$, we show the reversings 
$x^{-1} y \revstoer uv^{-1}$.
Since the set of
augmented relations has exactly one relation of the form $xu=yv$ for every pair
$\{x,y\}$ of generators, $x^{-1}y$ will always reverse, in one step, and uniquely, 
to a word of the form $uv^{-1}$ where $u$ and $v$ are positive words over the generators.
Furthermore, $x^{-1} y \revsto u v^{-1}$ precisely when $y^{-1} x\revsto v u^{-1}$,
so in the table we only include one of each of these. 

\begin{table}
  \centering  
\begin{enumerate}
\item{For any $p$, and any distinct $i$ and $j$, we have
$\left( a_p^{(i)} \right)^{-1} a_p^{(j)} \onerevsto
 a_p^{(i-1)} \left( a_p^{(j-1)} \right)^{-1}$};
\item{For any $p<q$, and any $i$ and $j$, we have
\begin{enumerate}
\item{$a_{pq}^{-1} a_{qp} \onerevsto  a_q^{(1)} a_q^{(0)} \left(a_p^{(1)} a_p^{(0)}\right)^{-1}$,}
\item{$a_{pq}^{-1} a_{q}^{(i)} \onerevsto a_q^{(i)} \left(a_p^{(i)}\right)^{-1} \quad
\mbox{ and } \quad
a_{qp}^{-1} a_{p}^{(i)} \onerevsto a_p^{(i)} \left(a_q^{(i-1)}\right)^{-1} $,}
\item{$a_{pq}^{-1} a_{p}^{(i)} \onerevsto a_q^{(i)} a_{pq}^{-1} \quad \quad 
\mbox{ and } \quad \quad \
a_{qp}^{-1} a_{q}^{(i)} \onerevsto a_p^{(i+1)} a_{qp}^{-1} $, \quad \quad and}
\item{$$\left( a_p^{(i)} \right)^{-1} a_q^{(j)} \onerevsto
\left\{\begin{array}{ll}
a_{pq} \left(a_p^{(i)}\right)^{-1} & \mbox{if $i=j$,} \\
a_p^{(i)} a_{qp}^{-1} & \mbox{ if $i = j+1$, } \\
a_p^{(i)} a_{pq} \left(a_q^{(j-1)} a_{qp}\right)^{-1} & \mbox{ if $i\neq j, j+1$ };
\end{array}\right.$$}
\end{enumerate}}
\item{For any $n$-cyclically ordered sequence $(p,q,r)$ and any $i$, we have
\begin{enumerate}
\item{$a_{pq}^{-1} a_{qr} \onerevsto  a_{qr} a_{pr}^{-1}, \quad
a_{pq}^{-1} a_{pr} \onerevsto  a_{qr} a_{pq}^{-1}$, \quad and \quad
$a_{qr}^{-1} a_{pr} \onerevsto  a_{pr} a_{pq}^{-1}$,}
\item{$a_{pq}^{-1} a_{r}^{(i)} \onerevsto a_{r}^{(i)} a_{pq}^{-1}$, }
\item{$a_{pr}^{-1} a_{qr} \onerevsto a_{pq} a_r^{(1)} a_r^{(0)} \left(a_{qr} a_p^{(1)} a_p^{(0)}\right)^{-1}$,}
\item{$$(a_{pr})^{-1} a_q^{(i)} \onerevsto
\left\{\begin{array}{ll}
a_{pq} a_r^{(i)} \left(a_{qr} a_p^{(i)}\right)^{-1} & \mbox{if $0 \leq p<q<r<n$,} \\
a_{pq} a_r^{(i)} \left(a_{qr} a_p^{(i-1)}\right)^{-1}& \mbox{ if $0 \leq q<r<p<n$, } \\
a_{pq} a_r^{(i+1)} \left(a_{qr} a_p^{(i)}\right)^{-1}& \mbox{if $0 \leq r<p<q<n$ };
\end{array}\right.$$}
\end{enumerate}}
\item{For any $n$-cyclically ordered sequence $(p,q,r,s)$, we have
\begin{enumerate}
  \item{$a_{pq}^{-1} a_{rs} \onerevsto  a_{rs} a_{pq}^{-1}$,}
  \item{$a_{ps}^{-1} a_{qr} \onerevsto  a_{qr} a_{ps}^{-1}$,}
  \item{$a_{pr}^{-1} a_{qs} \onerevsto  a_{pq} a_{qr} (a_{ps} a_{qr})^{-1}$,\ \ and}
  \item{$a_{ps}^{-1} a_{rq} \onerevsto  a_{pq} a_{qr} a_s^{(1)} a_s^{(0)}(a_{rp} a_{rs} a_q^{(1)} a_q^{(0)})^{-1}$.}
\end{enumerate}}
\end{enumerate}
\caption{
All the possible one-step reversings with respect to $\R_{e,n+1}$ for pairs of 
generators from $\A_{e,n+1}$. 
Exponents are always modulo $e$.
}
\label{onesteps}
\end{table}

Observe that the set of subscripts appearing on the right hand side of any of the
one-step reversings in Table~\ref{onesteps} is the same as the set of subscripts for
the left hand side. Furthermore, every ``exponent" (the $i$ appearing in $a^{(i)}$) on the
right hand side of a one-step reversing differs by at most one from an exponent
appearing on the left hand side; except for the exponents 0 and 1 which may appear
from nowhere.

The \emph{reversing calculation of $w$}, denoted $\rho_w$, where $w$ is a word
over $\Aer \cup \Aer^{-1}$, is defined to be the sequence 
of words obtained from one-step reversings, where at each step the left-most
occurrence of an (inverse generator, positive generator) pair is reversed,
(choosing this convention because it coincides with how we programmed the
reversing procedure in \gap), and continuing until no more reversings are 
possible, or otherwise continuing indefinitely. 
The last word in the sequence $\rho_w$, if it exists, is called the result of
the reversing calculation.
All the reversing calculations we undertook for the results we needed in this 
article terminated.
As remarked above, whenever the sequence $\rho_w$ is finite 
(in other words when the reversing procedure terminates) 
the result of $\rho_w$ is a word of the form $u v^{-1}$ 
where $u$ and $v$ are positive words.

\subsubsection*{Block mappings}

The \emph{exponent set} $I_w$  of a word is defined to be the set of exponents
appearing in the word; more precisely, if $w=a_{i_1} \cdots a_{i_f}$ 
where $a_{i_j} \in \Aer \cup \Aer^{-1}$, then
$$I_w = \bigcup_{j=1}^{f} I_{a_{i_j}}, 
\quad \mbox{ where } \quad 
I_{a_{pq}} := \varnothing, \quad I_{a_p^{(i)}} := \{ i \} \quad \mbox{ and } 
\quad I_{a^{-1}} := I_a \mbox{ for all $a \in \A$.} 
$$
The \emph{exponent set of a reversing calculation} is then the union of the
exponent sets for each word in the reversing calculation. By definition, the exponent
set will always be a subset of $\{0,1,\ldots,e-1\}$. The \emph{exponent blocks}
of an exponent set are the equivalence classes generated by the relation 
``$i \sim j$ if $|i-j|=1$ (modulo $e$)''. For example, given an exponent set
$I_{\rho_w} = \{0,1,2,3,6,7,11,12,13,17\}$ where $e=18$, the exponent blocks are
$\{17,0,1,2,3\}$, $\{6,7\}$ and $\{11,12,13\}$.
We will write this \mbox{$I_{\rho_w} = [17,3]_{18} \cup [6,7]_{18} \cup [11,13]_{18}$}, 
where $[i,k]_e$ denotes the $e$-cyclically ordered interval from $i$ to $k$. 
We will write $[e]$ to denote $\{0,\ldots,e-1\}$ with $e$-cyclic ordering.



Let $w$ be a word over $\Aer \cup \Aer^{-1}$, with reversing calculation $\rho_w$, 
and corresponding exponents $I_{\rho_w} = I_1 \cup \cdots \cup I_f$, in exponent
block decomposition. 
We have $I_l = [i_l, i_l+p_l-1]_e$ where $p_l$ is the size of $I_l$, and
$i_l$ its cyclically ordered least element. 
A \emph{block mapping} of $I_{\rho_w}$ to $[d]$ maps blocks to blocks,
acting identically on 0 and 1 if they appear in $I_{\rho_w}$.
More precisely, fix ${\mathbf j} = (j_1,\ldots,j_l)$ where each $j_t \in [d]$.
Define $\alpha_{\mathbf j}: I_{\rho_w} \rightarrow [d]$ as follows:
every $x \in I_{\rho_w}$ is of the form $x \equiv i_l + p$ modulo $e$ for some $l$ where
$0 \leq p < p_l$; then we set
$$\alpha_{\mathbf j}(x) :\equiv j_l + p \ \mod d \quad 
\mbox{ for } \quad x \equiv i_l + p \ \mod e \quad 
\mbox{ and } 0 \leq p < p_l.$$
Then $\alpha_{\mathbf j}$ is a \emph{block mapping} if
whenever $\{0,1\} \subseteq I_{\rho_w}$, we have 
$\alpha_{\mathbf j}(0) \equiv 0 \ \mod d$ and 
$\alpha_{\mathbf j}(1) \equiv 1 \ \mod d$.
Thus in particular:

\begin{lemma}
\label{betaplusk}
For $\alpha$ a block mapping: $I_{\rho_w} \rightarrow [d]$ and $k$ a positive integer,
whenever $[i,i+k] \subseteq I_l$ for some block $I_l$ of $I_{\rho_w}$,
we have \emph{$$\alpha(i+k) \equiv \alpha(i) + k \ \mod d.$$}
\end{lemma}

For example, given an exponent set $I_{\rho_w} = [17,3]_{18} \cup [6,7]_{18} \cup [11,13]_{18}$,
ordering the exponents from left to right $I_1, I_2, I_3$, we have that any 
block map $\alpha_{\mathbf j} : I_{\rho_w} \rightarrow [d]$ must have
$j_1 \equiv d-1 \ \mod d$, but that $j_2$ and $j_3$ can be chosen arbitrarily.
Thus $d=3$ and ${\mathbf j} = (2,1,1)$ gives rise to the block mapping
$$\begin{array}{rl}
\mbox{ on } I_1: \quad &
17 \bjmt 2 \quad \quad 0 \bjmt 0 \quad \quad 1 \bjmt 1
\quad \quad 2 \bjmt 2 \quad \quad 3 \bjmt 0, \\
\mbox{ on } I_2: \quad &
\ \, 6 \bjmt 1 \quad \quad 7 \bjmt 2, \\
\mbox{ and  on } I_3: \quad &
11 \bjmt 1 \quad \quad 12 \bjmt 2 \quad \quad 13 \bjmt 0.
\end{array}$$

We extend a block mapping $\alpha$ to the generators appearing in the reversing
calculation for $w$, and consequently homomorphically to words over these
generators, and so to the entire reversing calculation $\rho_w$, as follows:
$$\alpha(a_{pq}) := a_{pq}, \quad \alpha(a_p^{(i)}) := a_p^{(\alpha(i))}
\quad \mbox{ and } \quad \alpha(a^{-1}) := \alpha(a)^{-1} 
\quad \mbox{ for every $a$ in \Aer.}$$ 
Write $\alpha(\rho_w)$ for the sequence of words obtained by applying $\alpha$ 
to the sequence of words in the reversing calculation $\rho_w$ for $w$.

\begin{proposition}
Let $\rho_w$ be the reversing calculation with respect to \Rer of the word $w$. 
Then a block map $\alpha: I_{\rho_w} \rightarrow [d]$ 
sends $\rho_w$ to the reversing calculation of $\alpha(w)$ with respect to \Rdr \  --- that is, 
$$\alpha(\rho_w) = \rho_{\alpha(w)}.$$
\end{proposition}

\proof
It suffices to show that for every one-step reversing
$a^{-1} b \onerevsto_{\!\!\scriptscriptstyle{(e,r)}} u v^{-1}$ appearing in $\rho_w$, 
we have that
$\alpha(a^{-1} b) \onerevsto_{\!\!\scriptscriptstyle{(d,r)}} \alpha(u v^{-1})$.
We remark that $a$ and $b$ appear in $\rho_w$ only if 
all the exponents appearing in $a,b, u$ and $v$ lie in $I_w$.
For the numbering, we refer to Table~\ref{onesteps} where all the one-step
reversings are enumerated.
For the one-step reversings 3(a) and 4(a), (b) and (c) there is nothing
to prove, as $\alpha$ acts trivially on the generators appearing, and the 
relations are the same in \Rdr.
For the reversing 3(b), which is
$a_{pq}^{-1} a_{s}^{(i)} \onerevsto a_{s}^{(i)} a_{pq}^{-1}$, we have
$$
\alpha \left( a_{pq}^{-1} a_s^{(i)} \right) 
\  \equiv \  a_{pq}^{-1} a_s^{(\alpha(i))} 
\ \ \onerevsto_{\!\!\scriptscriptstyle{(d,r)}} \ \ a_s^{(\alpha(i))} a_{pq}^{-1} 
\  \equiv \ \alpha \left( a_s^{(i)} a_{pq}^{-1} \right).
$$
Let's now consider the one-step reversing numbered (1) in Table~\ref{onesteps}:
for any $p$, and any distinct $i$ and $j$, we have
$\left( a_p^{(i)} \right)^{-1} a_p^{(j)} \onerevsto
 a_p^{(i-1)} \left( a_p^{(j-1)} \right)^{-1}$. For this to occur in $\rho_w$
means that $i-1$ and $i$ (resp. $j-1$ and $j$) lie in some block $I_l$
(resp. $I_k$) of $I_w$. Then we have
$$\begin{array}{rlll}
\alpha \left( \left( a_p^{(i)} \right)^{-1} a_p^{(j)} \right) 
& \equiv & \left( a_p^{(\alpha(i))}  \right)^{-1} a_p^{(\alpha(j))} \\
& \onerevsto_{\!\!\scriptscriptstyle{(d,r)}} 
                 & a_p^{(\alpha(i)-1)} \left( a_p^{(\alpha(j)-1)} \right)^{-1} \quad
                 &\mbox{(where exponents are \mod $d$)}   \\
& \equiv & a_p^{(\alpha(i-1))} \left( a_p^{(\alpha(j-1))} \right)^{-1} 
                 &\mbox{(Lemma~\ref{betaplusk})}   \\
& \equiv & \alpha \left( a_p^{(i-1)} \left( a_p^{(j-1)} \right)^{-1} \right) .\\
\end{array}
$$
(The blocks $I_l$ and $I_k$ need not be distinct.)
Similar arguments go through for all the one-step reversings
2(b), (c) and (d),  and 3(d). 

The reversings remaining are 2(a), 3(c) and 4(d) which involve 
explicitly the exponents 0 and 1, and only these. By the definition of
the block mapping, since both 0 and 1 lie in $I_w$, we have that
$\alpha$ acts as the identity on these two exponents; 
hence its extension to words acts identically on 
the entirety of these relations.
$\Box$

\subsubsection{Results using \gap, applied to block mappings.}

In the following result, we collect together some results obtained 
using the \gap\ code. Recall that since all pairs of
generators have a corresponding one-step reversing, then a terminating
reversing calculation $\rho_w$ for any word $w$ results in a word
of the form $u v^{-1}$ where $u$ and $v$ are positive (but possibly
empty) words.

\begin{proposition}
\label{gapcalcs}
\begin{enumerate}
\item Let $p, q, r$ be any elements of $\{0,1,2\}$ (repeats allowed),
and denote 
$x = a_p^{(2)}$, $y = a_q^{(8)}$, $z = a_r^{(5)}$ 
and $w = x^{-1} z z^{-1} y$. Then the reversing calculation $\rho_w$ terminates, and
writing $u v^{-1}$ for the result of $\rho_w$, 
where $u$ and $v$ are positive words, we have
\begin{enumerate}
\item $w' :\equiv (xu)^{-1} yv\curvearrowright_{9,4} \ \varepsilon$, the empty word, and
\item  $I_{\rho_w}, I_{\rho_{w'}} \subseteq \{1,2\} \cup \{4,5\} \cup \{7,8\}$.
\end{enumerate}
\item Let $p,q,r,s$ be any elements of $\{0,1,2,3\}$ (repeats allowed),
let $\{x,y,z\} = \{a_{pq},a_r^{(2)},a_s^{(6)}\}$ (any ordering, but not repeating),
and let $w = x^{-1} z z^{-1} y$. Then  $\rho_w$ terminates, 
and writing $u v^{-1}$ for the result of $\rho_w$, where
$u$ and $v$ are positive words, we have
\begin{enumerate}
\item $w' :\equiv (xu)^{-1} yv \curvearrowright_{8,5} \ \varepsilon$, and
\item  $I_{\rho_w}, I_{\rho_{w'}} \subseteq \{1,2,3\} \cup \{5,6,7\}$.
\end{enumerate}
\item Let $p,q,r,s,t$ be any elements of $\{0,1,2,3,4\}$ (repeats allowed),
let $\{x,y,z\} = \{a_{pq},a_{rs},a_t^{(4)}\}$ (any ordering, but not repeating),
and let $w = x^{-1} z z^{-1} y$. Then $\rho_w$ terminates, and writing $u v^{-1}$ for
the result of $\rho_w$, where
$u$ and $v$ are positive words, we have
\begin{enumerate}
\item $w' :\equiv (xu)^{-1} yv \curvearrowright_{7,6} \ \varepsilon$, and
\item  $I_{\rho_w}, I_{\rho_{w'}} \subseteq \{0,1\} \cup \{3,4,5\}$.
\end{enumerate}
\item Let $p,q,r,s,t,f$ be any elements of $\{0,1,2,3,4,5\}$ (repeats allowed),
let $x=a_{pq}$, $y=a_{rs}$ and $z=a_{tf}$, 
and let $w = x^{-1} z z^{-1} y$. Then $\rho_w$ terminates, and writing $u v^{-1}$ for
the result of $\rho_w$, where
$u$ and $v$ are positive words, we have
\begin{enumerate}
\item $w' :\equiv (xu)^{-1} yv \curvearrowright_{4,7} \ \varepsilon$, and
\item  $I_{\rho_w}, I_{\rho_{w'}} \subseteq \{0,1,2\}$.
\end{enumerate}
\end{enumerate}
\end{proposition}

Recall that the completion condition with respect to \R\ holds for a triple
$x,y,z$ if whenever the word $x^{-1} z z^{-1} y$ reverses modulo \R\ 
to a word of the form $u v^{-1}$ where $u$ and $v$ are positive words, 
then $(xu)^{-1} yv \revstorelR \varepsilon$. 
Via well-chosen block mappings, we use the previous result to  
show that the completion condition holds for certain $r$ with no
restrictions on $e$.

\goodbreak

\begin{proposition}
\label{gencalcs}
\begin{enumerate}
\item For any $d$, the completion condition with respect to ${\cal R}_{d,4}$ 
holds for every triple 
$(a_p^{(i)},a_q^{(j)},a_r^{(k)})$ where $p,q,r \in \{0,1,2\}$ and 
$i,j,k \in [d]$.
\item For any $d$, the completion condition with respect to ${\cal R}_{d,5}$ 
holds for every triple consisting of
$a_{pq}$, $a_r^{(i)}$ and $a_s^{(j)}$ (any ordering) 
where $p,q,r,s \in \{0,1,2,3\}$ and 
$i,j,k \in [d]$.
\item For any $d$, the completion condition with respect to ${\cal R}_{d,6}$ 
holds for every triple consisting of 
$a_{pq}$, $a_{rs}$ and $a_t^{(i)}$ (any ordering)
where $p,q,r,s,t \in \{0,1,2,3,4\}$ and 
$i,j,k \in [d]$.
\item For any $d$, the completion condition with respect to ${\cal R}_{d,7}$ 
holds for every triple 
$(a_{pq}$, $a_{rs}$, $a_{tf})$ where $p,q,r,s,t,f \in \{0,1,2,3,4,5\}$ and 
$i,j,k \in [d]$.
\end{enumerate}
\end{proposition}

\proof
1. We continue to use the notation $w$, $u$ and $v$ for the words in  
Proposition~\ref{gapcalcs}, and define 
$$\overline{w} :\equiv \left( a_p^{(i)} \right)^{-1} a_r^{(k)} 
      \left( a_r^{(k)} \right)^{-1} a_q^{(j)}.$$ 
Define a block map 
$\bt: I_{\rho_w} \rightarrow [d]$ by
${\mathbf t} := (i-1,j-1,k-1)$; that is
$$\begin{array}{ll}
1\stackrel{\bt}{\longmapsto} i-1, \quad \quad \quad \quad &  2\stackrel{\bt}{\longmapsto} i, \\
4\stackrel{\bt}{\longmapsto} j-1, \quad \quad \quad \quad &  5\stackrel{\bt}{\longmapsto} j, \mbox{ and}\\ 
7\stackrel{\bt}{\longmapsto} k-1  \quad \quad \mbox{ and} & 8\stackrel{\bt}{\longmapsto} k.
\end{array}$$ 
Thus $\overline{w} \equiv \alpha(w)$.
Since $w \revsto_{\!\!\!\scriptscriptstyle{9,4}}\, u v ^{-1}$, we have
$\alpha(w) \revsto_{\!\!\!\scriptscriptstyle{d,4}}\, \alpha(u v^{-1}) 
\equiv \alpha(u) \alpha(v)^{-1}$. Furthermore, 
$$\alpha \left( \left(xu\right)^{-1} yv \right)  \equiv \left(a_p^{(i)} \alpha(u) \right)^{-1} a_q^{(j)} \alpha(v),$$
so $\left(xu\right)^{-1} yv \revsto_{\!\!\!\scriptscriptstyle{9,4}}\, \varepsilon$ implies
$$\left(a_p^{(i)} \alpha(u) \right)^{-1} a_q^{(j)} \alpha(v)
\revsto_{\!\!\!\scriptscriptstyle{d,4}}\, \alpha(\varepsilon) \equiv \varepsilon,$$
as required.
The arguments for the other three cases are the same as for 1., 
but with respect to different block mappings, being, respectively: 

\noindent
2. $\bt: \{1,2,3\} \cup \{5,6,7\} \rightarrow [d]$ where ${\mathbf t} = (i-1, k-1)$; that is
$$\begin{array}{lll}
1\stackrel{\bt}{\longmapsto} i-1, \quad \quad \quad \quad & 2\stackrel{\bt}{\longmapsto} i,
\quad \quad \quad \quad & 3\stackrel{\bt}{\longmapsto} i+1, \\
5\stackrel{\bt}{\longmapsto} k-1, &
6\stackrel{\bt}{\longmapsto} k, \quad \quad \mbox{ and } &
7\stackrel{\bt}{\longmapsto} k+1.\\
\end{array}$$ 
3. $\bt: \{0,1\} \cup \{3,4,5\} \rightarrow [d]$ where ${\mathbf t} = (0,i-1)$; that is
$$\begin{array}{lll}
0\stackrel{\bt}{\longmapsto} 0, \quad \quad \quad \quad 
& 1\stackrel{\bt}{\longmapsto} 1,\quad \quad \quad \quad& \\
3\stackrel{\bt}{\longmapsto} i-1,  \quad \quad \quad \quad &
4\stackrel{\bt}{\longmapsto} i, \quad \quad \mbox{ and } &
5\stackrel{\bt}{\longmapsto} i+1.\\
\end{array}$$ 
4. $\alpha_4: \{0,1,2\} \rightarrow [d]$ where  ${\mathbf t} = (0)$; that is
$\quad
0\stackrel{\bt}{\longmapsto} 0, \quad 
1\stackrel{\bt}{\longmapsto} 1 \quad \mbox{ and } \quad 
2\stackrel{\bt}{\longmapsto} 2.
\quad \quad \quad$
$\Box$

It remains to extend these results to the case of general $r$; having done
which we will have shown that each presentation $\langle \Aer \mid \Rer \rangle$
is right complete.

\subsubsection*{Principal subscript maps}

Let $S_w$ denote the set of subscripts appearing in the word $w$; that is
$$S_{a_pq} = \{p,q\}, \quad S_{a_p^{(i)}} = \{p\} \quad \mbox{ and } \quad
S_{a^{-1}} = S_a \quad \mbox{ for all } a \in \A,$$
extending to words in the obvious way: $S_{wx} = S_w \cup S_x$.
By inspection of Table~\ref{onesteps}, if $x ^{-1} y \onerevsto u v^{-1}$
then $S_{x^{-1}y} = S_{uv^{-1}}$. Thus the subscript set of a reversing calculation 
is exactly the subscript set of its first word.

Define the \emph{principal subscript map} $\sigma_w$ of a word $w$ as follows: 
let $S_w = \{p_0,\ldots,p_f\}$ where $0 \leq p_0 < \cdots < p_f$, and define
$\sigma_w (p_i) = i$. Thus any ordering $(p_{i_1},\ldots,p_{i_k})$ of a subset  of $S_w$ 
is cyclically ordered precisely when $({i_1},\ldots,{i_k})$ is cyclically ordered.
As an example, the word $w = a_{1,7}^{-1} a_{6,9}^{~} a_{6,9}^{-1} a_{4}^{(15)} $ has
$S_w = \{0,4,6,7,9\}$ and corresponding principal subscript map
$$\sigma_w : \quad 1 \mapsto 0, \quad 4 \mapsto 1, \quad 6 \mapsto 2, \quad 7 \mapsto 3 \quad
\mbox{ and } \quad 9 \mapsto 4. $$ 
We then extend this in the obvious way to arbitrary words whose subscript set is 
a subset of $S_w$: 
$$\sigma_w(a_{p_s p_t}) := a_{st}, \quad \sigma_w (a_{p_s}^{(i)}) := a_{s}^{(i)} \quad
\mbox{ and } \quad \sigma_w(a^{-1}) := \sigma(a) \quad
\mbox{ for every } a \in \Aer.$$
For the example above, 
this gives $\sigma_w(w) = a_{0,3}^{-1} a_{2,4}^{~} a_{2,4}^{-1} a_{1}^{(15)} $.
In general, for $w'$ a word whose subscripts lie inside $S_w$,
$\sigma_w(w')$ is a word whose subscripts are all strictly less than $|S_w|$.

Because the relations \Rer\ are defined in terms of cyclic order on the subscripts,
inspection of Table~\ref{onesteps} will convince the reader that
for principal subscript map $\sigma_w$ whose domain contains $S_{x^{-1} y}$,
$$ x^{-1} y \ \onerevsto_{\scriptscriptstyle{e,r}} \ u v^{-1} \ \Leftrightarrow \ 
\sigma(x^{-1} y) \ \onerevsto_{\scriptscriptstyle{e,r'}} \ \sigma(u v^{-1})$$
where $r' = |S_w|$.

\begin{proposition}
\label{allr}
Let $r' = | S_{x y z} |$ for $x,y,z \in \Aer$. Then
the completion condition holds for $(x,y,z)$ modulo \Rer\ if and only if
the completion condition holds for $(\sigma(x), \sigma(y), \sigma(z))$, 
which is a triple over ${\cal A}_{e,r'}$, modulo ${\cal R}_{e,r'}$, 
where $\sigma$ is the principal subscript map on $S_{xyz}$.
\end{proposition}

Observe that in Proposition~\ref{gencalcs} we proved
that completion holds for each different type of triple with respect
to the corresponding principal subscript set. 
Combining Propositions~\ref{gencalcs} and~\ref{allr}, we obtain:

\begin{theorem}
\label{complete}
The augmented presentation $\langle \A_{e,r} \mid \R_{e,r} \rangle$ 
is right complete for any choice of $e$ and $r$. 
Thus the $B^+$ is left cancellative, and the partial order
defined by left division is a lattice.
\qed
\end{theorem}

Invoking Theorem~\ref{revisom} we observe that \Bprev\ is also left
cancellative and its left division poset is a lattice.
Applying the anti-isomorphism \rev, we deduce that $B^+$ must be right cancellative, and 
that the partial order defined by right division must be a lattice. 
Hence we have:

\begin{theorem}
\label{leftandright}
The monoid $B^+$ is cancellative and the posets
defined by left and right division are lattices. \qed
\end{theorem}

\subsubsection*{Left completeness}

A presentation is left complete if its reverse is right complete, where
the reverse presentation is the presentation with the same generators
and all relations of the form $\rev(\rho_1) = \rev(\rho_2)$ whenever
$\rho_1 = \rho_2$ is a relation from the original presentation. 

We observe that the presentation which we showed to be left complemented and left complete
is in general neither right complemented (since it contains, for example, relations of type $\C_2$
which are of the form $ux = vx$) nor right complete (since it contains, for example, 
no relation of the form $u a_{pr} = v a_{q}^{(i)}$, for $p<q<r$, while 
by Theorem~\ref{leftandright}, $a_{pr}$ and $a_{q}^{(i)}$ do have a left common multiple, 
excluding completeness).
We suspect that there is no presentation for $B^+$ on our generating set which is simultaneously
left and right complete. Our point here is to show that
$\langle \A \mid \rev(\chi(\R^+)) \rangle$ is a presentation 
for $B^+$ which is left complete. This result is not necessary in the sequel, but is
included for reasons of completeness.

Recall that we defined a permutation $\chi$ of \A\ on page~\pageref{chidef}.
Since the augmenting relations of types $\C_1, \ldots, \C_6$ 
are consequences of the relations \R, the relations
of types $\chi(\C_1),\ldots,\chi(\C_6)$ are consequences
of the relations $\chi(\R) = \rev(\R)$. Thus
$$\langle \A \mid \chi(\R^+) \rangle = \langle \A \mid \chi(\R) \rangle = \Bprev.$$ 
This presentation is left complemented, because
$\langle \A \mid \chi(\R^+) \rangle$ is left complemented.
Furthermore, since 
$$\begin{array}{rcccl}
a^{-1}b \onerevstorelR uv^{-1} &\Longleftrightarrow& (au = bv) \in \R \\
&\Longleftrightarrow& (\chi(au) = \chi(bv) ) \in \chi(\R) 
&\Longleftrightarrow& \chi(a^{-1}b) \onerevstorelchiR \chi(uv^{-1})
\end{array}
$$
we have that
\begin{lemma}
\label{chipreservesreversing}
The operator $\chi$ on $\A^*$ preserves word reversing: that is, 
for any word $w$ over \A, $$\chi(\rho_w) = \rho_{\chi(w)}.$$ 
\end{lemma}
Hence application of $\chi$ maps the verification of the completeness condition for 
$\langle \A \mid \R^+ \rangle$ to the verification of the completeness condition for
$\langle \A \mid \chi(\R^+) \rangle$. Thus $\langle \A \mid \chi(\R^+) \rangle$
is right complete and by applying $\rev$ we have that 
$\langle \A \mid \rev(\chi(\R^+)) \rangle$ is left complete.
Thus in summary,

\begin{theorem}
Let \emph{$\Q^+:=\rev(\chi(\R^+))$}. Then the presentation 
$\langle \A \mid \Q^+ \rangle$
is a left complemented and left complete presentation for $B^+$. \qed
\end{theorem}


\subsection{Garside element}
\label{sec:garelement}

To show that the monoid $B^+$ is Garside, it remains to show that it has
a Garside element -- an element with the property that its set of left divisors
is equal to its set of right divisors. We will use the following general result.

\begin{lemma}
\label{generalemma}
Let $M$ be a cancellative monoid and suppose that 
$\overline{\, \White{\cdot} \, } : M \longrightarrow M$ is an automorphism of $M$ 
for which there is a fixed element $\Delta$ of $M$ such that
$m \Delta = \Delta \overline{m}$ for all $m$ in $M$.
Then the set of left and right divisors of $\Delta$ are equal.
\end{lemma}

\begin{proof}
Suppose that $\Delta = mn$. Then 
$mn \overline{m} = \Delta \overline{m} = m \Delta$.
Left cancelling $m$ we have that 
$n \overline{m} = \Delta$, showing that every right divisor
of $\Delta$ is also a left divisor. Furthermore, since 
$\overline{\Delta} = \Delta$, applying the inverse of
the automorphism to $n \overline{m} = \Delta$ shows that
every left divisor is also a right divisor.
\qed
\end{proof}

We will use the notation $\binom{q}{p}$ to represent the 
product $a_{p, p+1} a_{p+1, p+2} \cdots a_{q-1,q}$, where the subscripts are 
taken modulo $n$; if $p=q$ then this is 
the empty word, representing the identity;
$\binom{p+1}{p} = a_{p,p+1}$; and if $p>q$ then 
$\binom{q}{p} = a_{p, p+1} \cdots a_{n-1,0} a_{0,1} \cdots a_{q-1,q}$.
Clearly if  $(p,q,r)$ is cyclically ordered then $\binom{q}{p} \binom{r}{q} = \binom{r}{p}$.
The following relations will be useful in the sequel.

\begin{lemma}
\label{binomrels}
\begin{enumerate}
\item For any $p \leq q$ and any $i$,  $a_p^{(i)} \binom{q}{p} = \binom{q}{p} a_q^{(i)}$,
\item For any $(p,q,r)$ cyclically ordered, 
$a_{p,r} \binom{q}{p} = \binom{q}{p} a_{q,r}$ and $a_{r,p} \binom{q}{p} = \binom{q}{p} a_{r,q}$.
\end{enumerate}
\end{lemma}

\begin{proof}
Recall that by definition,
$\binom{p+1}{p} = a_{p, p+1}$.  
Translating  relations $(\R_5)$ and $(\R_3)$ respectively into this language, we have
\begin{enumerate}
\item for all $p$, $a_p^{(i)} \binom{p+1}{p} = \binom{p+1}{p} a_{p+1}^{(i)}$, and
\item for $(p,p+1,r)$ cyclically ordered, $a_{p,r} \binom{p+1}{p} = \binom{p+1}{p} a_{p+1,r}$ and
$a_{r,p} \binom{p+1}{p} = \binom{p+1}{p} a_{r,p+1}$.
\end{enumerate}
Since each (non-empty) $\binom{q}{p}$ is of the form $\prod_{t=p}^{q-1} \binom{t+1}{t}$,
the result follows immediately.
\qed
\end{proof}

Let $\alpha := a_0^{(1)} a_0^{(0)} \binom{n-1}{0}$ (this is exactly the image of 
the element $\beta$ of~\cite{bmr} under our isomorphism of the next section).
In this section, we will show that this is the lcm of the generators \A, 
and is a Garside element for the braid monoid.

We begin by observing that the permutation $\rho$ of \A\  given by
$${a_{pq}} \stackrel{\rho}{\longmapsto} a_{p+1,q+1} \quad \mbox{ and } \quad
{a_p^{(i)}} \stackrel{\rho}{\longmapsto} \left\{ 
\begin{array}{ll}
a_{0}^{(i+1)} &\mbox{ if $p = n-1$, and }\\
a_{p+1}^{(i)} &\mbox{ otherwise, }
\end{array} \right.
$$  
(subscripts considered modulo $n$, exponents modulo $e$)
extends to an automorphism of the braid monoid. It is enough to check 
that $\rho(\R) =  \R$; by finiteness of \R\ it suffices to show
$\rho(\R) \subseteq \R$, which is clear by inspection.
The action of $\rho$ on \A\ has $n$ orbits, with representatives
$a_{0,q}$ for each $q$ in $0<q<n$, 
together with $a_0^{(0)}$, where  
$$\rho^p(a_{0,q}) = a_{p,p+q}, \quad \mbox{ and for $0\leq p<n$, }  \quad
\rho^{ni+p}(a_{0}^{(0)}) = a_{p}^{(i)}.$$

Clearly, for any $p,q$, we have $\binom{q}{p} \stackrel{\rho}{\mapsto} \binom{q+1}{p+1}$. 
Also, writing $t_q$ for $a_q^{(1)} a_q^{(0)}$, observe that $\rho(t_q) = t_{q+1}$.
Using these facts, along with Lemma~\ref{binomrels}.2 and $\C_6$, we see that the element 
$\alpha$ is preserved by $\rho$:
$$\begin{array}{c}
\rho(\alpha) = \rho\left(t_0 \binom{n-1}{0} \right) =  t_{n-1} \binom{n-1}{n-2}
 =  t_{n-1} a_{n-1,0} \binom{n-2}{0}  = t_0 a_{0,n-1 }\binom{n-2}{n-1} 
 =  t_0 \binom{n-2}{0} a_{n-2, n-1}  = \alpha.
\end{array}
$$
\begin{proposition}
\label{abeta}
There is a permutation $\overline{\cdot}$ of \A, such that for all $a$ in \A, 
$$\alpha = a \alpha_a = \alpha_a \overline{a}$$
for some element $\alpha_a$. 
\end{proposition}

\proof
Define $\gamma_0 = a_0^{(e-1)} \binom{n-1}{0}$ and 
$\gamma_q = t_q \binom{q-1}{0}  \binom{n-1}{q}$ for $0<q<n$. Then we have
$$\begin{array}{rcl}
\alpha &= & a_0^{(0)} \gamma_0  \\
& = & a_0^{(e-1)} a_0^{(e-2)} \binom{n-1}{0} 
=   a_0^{(e-1)} \binom{n-1}{0} a_{n-1}^{(e-2)} = \gamma_0 a_{n-1}^{(e-2)}, \mbox{ and } \\
~\\
\alpha & = & t_0 \binom{q-1}{0} a_{q-1,q} \binom{n-1}{q} \\
& = & t_0 a_{0q} \binom{q-1}{0} \binom{n-1}{q} \ = \  a_{0q} t_q  \binom{q-1}{0} \binom{n-1}{q} 
\ = \ a_{0q} \gamma_q \\
& = & t_q a_{q0} \binom{q-1}{0} \binom{n-1}{q} \ = \  t_q  \binom{q-1}{0} a_{q,q-1} \binom{n-1}{q} 
\ = \  t_q  \binom{q-1}{0} \binom{n-1}{q} a_{n-1,q-1}  \ = \  \gamma_q a_{n-1,q-1}.
\end{array}$$
That is, for each of the $\rho$-orbit representatives $x$, $\alpha = x \alpha_x = \alpha_x \rho^k(x)$
where $k = -(n+1)$. Applying $\rho$ an appropriate number of times, and recalling that 
$\rho(\alpha)= \alpha$, gives the equation
$$\alpha = x \alpha_x = \alpha_x \rho^{-(n+1)}(x)$$
for each generator $x$, for some $\alpha_x$. Defining $\overline{x} := \rho^{-(n+1)}(x)$,
we obtain the stated result. 
\qed

This proposition also shows that $\alpha$ is a common multiple of \A.
Since the braid monoid has a complete presentation, every set with a common
multiple has a least common multiple, so in particular \A has a least common multiple.

\begin{proposition}
\label{betaislcm}
The element $\alpha$ is the least common multiple of the generators.
\end{proposition}

\begin{proof}
We use the fact that if a generator $a$ does not divide $x$ but $a$ divides $xb$
for some generator $b$ then $xb$ is the least common multiple of $x$ and $a$.
This follows directly from the fact that the relations are homogeneous.

We show firstly that $\binom{r-1}{0}$ is the least common multiple of 
$\{a_{pq} | 0 \leq p < q < r\}$: this proceeds via induction  on $r$.
The base case $r=2$ holds since $a_{01}$ is the lcm of $\{a_{01}\}$.
For any $r$, it is easily seen that $\binom{r}{0}$ is 
a common multiple of $\{a_{pq} | 0 \leq p < q < r+1\}$. Furthermore,
$\binom{r}{0} = \binom{r-1}{0} a_{r-1, r}$. Since the subscript $r$ does
not appear anywhere in $\binom{r-1}{0}$, $a_{r-1,r}$ does not divide
$\binom{r-1}{0}$ (observed in Section~\ref{sec:pres}); but $a_{r-1,r}$
does divide $\binom{r}{0}$, so by the observation of the previous paragraph
we have that $\binom{r}{0}$ is the lcm of $\binom{r-1}{0}$ and $a_{r-1,r}$.
But being a common multiple of $\{a_{pq} | 0 \leq p < q < r+1\}$
implies that it is in fact the least common multiple of the whole set.

Thus in particular, $x := \binom{n-1}{0}$ is the least common multiple of 
$\{a_{pq} | 0 \leq p < q < n\}$. By the same argument, since 
$a_0^{(1)}$ does not divide $x$ (the only relations which 
could be applied being $(\R_1), (\R_2), (\R_3)$ which contain no $a_p^{(i)}$) 
but $x a_{n-1}^{(1)} = a_0^{(1)} x$, we have that
$a_0^{(1)} x$ is the lcm of $a_0^{(1)}$ and $\{a_{pq} | 0 \leq p < q < n\}$.

Finally, since the lcm of $a_0^{(1)}$ and $a_0^{(0)}$ is 
$a_0^{(1)} a_0^{(0)} = a_0^{(0)} a_0^{(e-1)}$, we have that
any common multiple of $a_0^{(1)}$ and $a_0^{(0)}$ must be divisible by 
$a_0^{(1)} a_0^{(0)}$; since $x$ is not divisible by $a_0^{(0)}$
then $a_0^{(1)} x$ is not either. But since 
$a_0^{(1)} x a_{n-1}^{(0)} = \alpha$ is divisible by all the 
generators, then it is divisible by $a_0^{(0)}$, and hence by the observation
of the first paragraph of this proof, is the lcm of $\{a_{pq} | 0 \leq p < q < n\}$, 
$a_0^{(1)}$ and $a_0^{(0)}$. Thus $\alpha$ is the lcm of all the generators.
\qed
\end{proof}

\begin{theorem}
\label{isgarside}
The least common multiple $\alpha$ of the generators \A\ is a Garside element for the braid monoid.
\end{theorem}

\proof
We have that $\alpha$ is the least common multiple of the generators and that
$\alpha = a \alpha_a = \alpha_a \overline{a}$ for each $a$. Thus 
$a \alpha = a \alpha_a \overline{a} = \alpha \overline{a}$, for each $a$, and we
can invoke Lemma~\ref{generalemma}. 
\qed


\section{The presentation \emph{does} present the braid group of $G(e,e,r)$}
\label{sec:pres}

In this section we show that the group defined by the presentation introduced here
is indeed the braid group of type $G(e,e,r)$.

\subsection*{The presentation of \cite{bmr}}

From~\cite{bmr} we know there is a presentation of the braid group of $G(e,e,n+1)$ as follows:

\noindent
\emph{Generators:} ${\cal T} = \{\tau_2, \tau'_2, \tau_3, \ldots, \tau_n, \tau_{n+1}\}$

\noindent
\emph{Relations:}  The commuting relations are: 
$\tau_i \tau_j = \tau_j \tau_i$ whenever $|i-j| \geq 2$,
together with $\tau'_2 \tau_j = \tau_j \tau'_2$ for all $j \geq 4$.
The others are:
$$\begin{array}{rcll}
\langle \tau_2 \tau'_2 \rangle^e &= &\langle \tau'_2 \tau_2 \rangle^e \\
\tau_i \tau_{i+1} \tau_i &= &\tau_{i+1} \tau_i \tau_{i+1} &\mbox{ for }i = 2, \ldots, n\\
\tau'_2 \tau_3 \tau'_2 &= &\tau_3 \tau'_2 \tau_3\\
\tau_3 \tau_2 \tau'_2 \tau_3 \tau_2 \tau'_2 &= &\tau_2 \tau'_2 \tau_3  \tau_2 \tau'_2 \tau_3
\end{array}$$
where the expression $\langle ab \rangle^k$ denotes the alternating product
$aba\cdots$ of length $k$.

\noindent
We call these the BMR-generators and the BMR-relations.
Recall that we started by looking for an answer to the following:

\emph{(1) Is the natural morphism $M \rightarrow B$ injective?}

\emph{(2) Do we have 
$$B = \{\ \alpha^n b \ | \ (n \in \mathbb{Z} ) \ ( b \in M )\ \}?$$}

\noindent where $M$ denotes the monoid defined by the above presentation.

For ease on the eye, we will write $1$ for $\tau_2'$, $2$ for $\tau_2$ and
$i$ for $\tau_i$ for $i=3,\ldots, n-1$.
Observe that $$
\begin{array}{rcl}
2\ 1  3213 \langle 21 \rangle^{e-2} 
        & = & 321321 \langle 21 \rangle^{e-2} \ = \ 3213 \langle 21 \rangle^{e} 
         \ = \ 3213 \langle 12 \rangle^{e} \\
         & = & 3213 \ 1 \langle 21 \rangle^{e-1} \ = \ 32 313 \langle 21 \rangle^{e-1} 
         \ =  \ 2 \ 3213 \langle 21 \rangle^{e-1} .
\end{array}
$$
However $1 3213 \langle 21 \rangle^{e-2}$ and 
$3213 \langle 21 \rangle^{e-1}$ have no subwords
appearing in ${\cal R}_0$, so are in singleton (hence distinct) ${\cal R}_0$-equivalence classes. Thus
the monoid $M$ is not cancellative, and so does not embed in $B$.

\begin{proposition}
The monoid defined by the presentation for the braid group $B$
of $G(e,e,r)$ given in~\cite{bmr} is not 
cancellative, and hence does not embed in the braid group. 
Denote by ${\cal S}$ the set of relations obtained by adjoining the relation 
$$1 3 2 1 3 \langle 2 1 \rangle^{e-2} 
= 3 2 1 3 \langle 2 1 \rangle^{e-1}$$
to the set of BMR relations for $B$. (Since this relation holds in $B$,
adjoining the relation does not change the group thus presented.)
Then  
$2 1^n 3 2 1 3 \langle 2 1 \rangle^{e-2} 
\Seq 3 2 1 3 2 ^n 1 \langle 2 1 \rangle^{e-2}$ for every $n$,
while $2 1^n 3 2 1 3 $ and 
$3 2 1 3 2^n 1$ are in distinct $\cal S$-classes for every $n>1$. 
Thus $M({\cal T}, {\cal S})$ is not cancellative either.
\end{proposition} 

\proof
The new relation can be written $1 w = w x $ where $w$ is $3213 \langle 21 \rangle^{e-2}$ 
and $x$ is $2$ if $e$ is even, and $1$ if $e$ is odd. For the same letter $x$,
$\langle 21 \rangle^e x \equiv \langle 21 \rangle^{e+1} \equiv 2 \langle 12 \rangle^e$
hence for any $n$ we have that 
$$21^n 3213\langle 21 \rangle^{e-2} \equiv 21^n w 
\Seq 21 w x^{n-1}  \Req 3213 \langle 21 \rangle^e x^{n-1}
\Seq 3213 2^{n-1} \langle 21 \rangle^e  \equiv  3213 2^n 1 \langle 21 \rangle^{e-2}.$$
However for every $n>1$, $21^n 3213$ and $3213 2^n 1$ are in singleton \R-equivalence classes; 
so the monoid $M({\cal T},{\cal S})$ is not cancellative, and so does not embed in $B$.
\qed

In the next subsection we show that the new presentation we have is a presentation
for the braid group of $G(e,e,r)$, giving a monoid answering affirmatively 
to all the desired properties. Using the facts that with the new presentation
we have a solution to the word problem and can calculate least common multiples, 
we show at the end of the section that there can be no finite presentation for
the braid group of $G(e,e,r)$ on the BMR-generators for which the 
corresponding monoid embeds in the group.

\subsection*{The new presentation presents $G(e,e,r)$}

Let $B$ denote the group defined by the BMR-presentation, and $G$ denote the group 
defined by our presentation, as given in Section~\ref{statement}. 
Define $\varphi: B \rightarrow G$ as follows: on generators, 
$$\varphi: \left\{
 \begin{array}{rcll} 
\tau_2 &\mapsto & a_0^{(1)} \\
\tau'_2 &\mapsto & a_0^{(0)} \\
\tau_i &\mapsto & a_{i-3, i-2} &\mbox{ for } i = 3,\ldots,n+1,
\end{array} \right.
$$
and extend homomorphically. To see this is well-defined it suffices to show that whenever
$\rho_1 = \rho_2$ is a defining relation in $B$, then $\varphi(\rho_1)$ can be transformed
into $\varphi(\rho_2)$ using the defining relations of $G$. This can be seen for the 
commuting relations simply by inspection -- they appear as incarnations of 
$(\R_1)$, $(\R_2)$  and $(\R_4)$.


Observe that 
$\langle a_0^{(1)} a_0^{(0)} \rangle^{e-1} 
= (a_0^{(e-1)} a_0^{(e-2)})  \ \cdots \ a_0^{(1)}$, so
$$\langle a_0^{(0)} a_0^{(1)} \rangle^e = 
a_0^{(0)} \langle a_0^{(1)} a_0^{(0)} \rangle^{e-1}
 = a_0^{(0)} (a_0^{(e-1)} a_0^{(e-2)})  \ \cdots \ a_0^{(1)}
= \langle a_0^{(1)} a_0^{(0)} \rangle^{e}.$$
Thus we have 
$\langle \tau_2 \tau'_2 \rangle^e 
\stackrel{\varphi}{\mapsto}  \langle a_0^{(1)} a_0^{(0)} \rangle^e 
= \langle a_0^{(0)} a_0^{(1)} \rangle^e 
\stackrel{\varphi}{\mapsfrom} \langle \tau'_2 \tau_2 \rangle^e$
as required.
 
From $(\R_3)$ we have, for $i \geq 3$, 
$$\tau_i \tau_{i+1} \tau_i \stackrel{\varphi}{\mapsto}a_{j-1,j} a_{j,j+1} a_{j-1,j}
\stackrel{{\cal R}_3}{=} a_{j,j+1} a_{j-1,j+1} a_{j-1,j}
\stackrel{{\cal R}_3}{=} a_{j,j+1}  a_{j-1,j} a_{j,j+1} 
\stackrel{\varphi}{\mapsfrom} \tau_{i+1} \tau_i \tau_{i+1},$$
where $j = i-3$. Using $(\R_5)$, we get
$$\tau_2 \tau_{3} \tau_2 \stackrel{\varphi}{\mapsto}a_0^{(1)} a_{01} a_0^{(1)}
\stackrel{{\cal R}_5}{=} a_{01} a_{1}^{(1)} a_0^{(1)}
\stackrel{{\cal R}_5}{=} a_{01} a_0^{(1)} a_{01}
\stackrel{\varphi}{\mapsfrom} \tau_{3} \tau_2 \tau_{3},$$
and similarly replacing $\tau_2$ with $\tau'_2$ and 
$a_0^{(1)}$ with $a_0^{(0)}$.
Finally, using the equation after $(\C_6)$,
$$\tau_3 \tau_2 \tau'_2 \tau_3 \tau_2 \tau'_2 
\stackrel{\varphi}{\mapsto} a_{01} t_0 a_{01} t_0
{=} t_0 a_{01} t_0 a_{01}
\stackrel{\varphi}{\mapsfrom} 
\tau_2 \tau'_2 \tau_3  \tau_2 \tau'_2 \tau_3.$$
Thus we have that $\varphi$ is well-defined.


For any $0 \leq p < q < n$, we use the notation 
$$\begin{array}{c}
\binom{q}{p} := \left\{
\begin{array}{ll}
\varepsilon \mbox{ (the empty word) } & \mbox{if }p = q, \\
\tau_{p+3} \tau_{p+4} \cdots \tau_{q+2} & \mbox{otherwise.}
\end{array} \right.
\end{array}
$$
Thus, for example, $\binom{n-1}{0} = \tau_3 \tau_4 \cdots \tau_{n+1}$.
Observe that $\binom{i+1}{i} = \tau_{i+3}$ for any $0 \leq i < n$, and
this maps to $a_{i, i+1}$ under $\varphi$. Thus  $\binom{q}{p}$ maps
under $\varphi$ to the element $\binom{q}{p}$ of $G$ 
defined in Section~\ref{sec:garelement}. 
Define $\beta := \tau_2 \tau_2' \binom{n-1}{0}$ (which is 
precisely the generator of the centre given in \cite{bmr});
thus $\beta \stackrel{\varphi}{\rightarrow} \alpha$, the 
least common multiple of the generators of the braid monoid.

In Section~\ref{sec:garelement} we observed that defining
$\overline{x} := \alpha^{-1} x \alpha$
restricts to a permutation on \A, in particular such that
$$\overline{a_{pq}} = a_{p-1,q-1} \quad  \mbox{ and } \quad 
\overline{a_p^{(i)}} =  \left\{
\begin{array}{ll}
a_{n-1}^{(i-2)} &\mbox{if $i=0$}, \\
a_{p-1}^{i-1} &\mbox{otherwise.}
\end{array}
\right.$$
We will mimic this in $B$.

For $0\leq p < q < n$ define $b_{pq} := \bpq$, $b_{qp} := \bqp$, 
and 
$$ b_0^{(i)} = 
\left\{\begin{array}{ll}
1 & \mbox{if $i = 0$,}\\
\boi & \mbox{for $i=1,\ldots e-1$,}
\end{array} \right.$$
where we again write $1$ for $\tau_2'$, $2$ for $\tau_2$.
Lastly, define $b_p^{(i)} := \bop^{-1} b_0^{(i)} \bop$.

\begin{lemma}
\label{alpharot}
$$\beta^{-1} b_{pq} \beta = b_{p-1,q-1} \quad  \mbox{ and } \quad 
\beta^{-1} b_p^{(i)} \beta =  \left\{
\begin{array}{ll}
b_{n-1}^{(i-2)} &\mbox{if $i=0$}, \\
b_{p-1}^{i-1} &\mbox{otherwise.}
\end{array}
\right.$$
\end{lemma}

\proof
\label{proofalpharot}
The proof is by calculation of the various cases.
\begin{enumerate}
\item We show $\beta^{-1} b_{pq} \beta = b_{p-1,q-1}$ for $p<q$;
there are two cases, in the second, $p>0$.
$$\begin{array}{lrcl}  
(i) &\beta b_{n-1, q-1} \beta^{-1} 
&= &21 \binom{n-1}{0} \binom{q-1}{0}^{-1} \binom{n-1}{0}^{-1} 1^{-1} 2^{-1} 321
\binom{n-1}{0} \binom{q-1}{0} \binom{n-1}{0}^{-1} 1^{-1} 2^{-1}\\
&&= &21 \binom{q}{1}^{-1} \binom{n-1}{0} \binom{n-1}{0}^{-1} 1^{-1} 2^{-1} 321
\binom{n-1}{0} \binom{n-1}{0}^{-1} \binom{q}{1} 1^{-1} 2^{-1}\\
&&= &\binom{q}{1}^{-1} 21 1^{-1} 2^{-1} 321 1^{-1} 2^{-1}\binom{q}{1}\\
&&= &\binom{q}{1}^{-1} 3\binom{q}{1} \ = \ \binom{q}{1}^{-1} \binom{q}{0}
\ = \ b_{0,q}, \\
\\
(ii) &\beta b_{p-1, q-1} \beta^{-1}
&= &21 \binom{n-1}{0} \binom{q-1}{p}^{-1} \binom{q-1}{p-1} \binom{n-1}{0}^{-1} 1^{-1} 2^{-1}\\
&&= &21 \binom{q}{p+1}^{-1} \binom{n-1}{0} \binom{n-1}{0}^{-1} \binom{q}{p} 1^{-1} 2^{-1}\\
&&= &\binom{q}{p+1}^{-1} 21 1^{-1} 2^{-1}\binom{q}{p} 
\ = \ \binom{q}{p+1}^{-1} \binom{q}{p}
\ = \ b_{p,q}. 
\end{array}$$
\item Next we show that $\beta^{-1} b_{qp} \beta = b_{q-1,p-1}$ for $p<q$;
again two cases, and in the second, $p>0$.
$$\begin{array}{lrcl}  
(i) &\beta b_{q-1, n-1} \beta^{-1} 
&= &21 \binom{n-1}{0} \binom{n-1}{q}^{-1} \binom{n-1}{q-1} \binom{n-1}{0}^{-1} 1^{-1} 2^{-1}\\
&&= &21 \binom{q}{0}  \binom{q-1}{0}^{-1} 1^{-1} 2^{-1}\\
&&= &21  \binom{q}{1}^{-1} \binom{q}{0} 1^{-1} 2^{-1}\\
&&= &  \binom{q}{1}^{-1} 21 3 1^{-1} 2^{-1} \binom{q}{1} \\
&&= &  \binom{q}{0}^{-1} 321 3 1^{-1} 2^{-1} 3^{-1} \binom{q}{0} \\
&&= &  \binom{q}{0}^{-1} 1^{-1} 2^{-1} 321 3 3^{-1} \binom{q}{0} 
\ = \ \binom{q}{0}^{-1} 1^{-1} 2^{-1} 321 \binom{q}{0} 
\ = \ b_{q0}\\
\\
(ii) &\beta b_{q-1, p-1} \beta^{-1} 
&= &21 \binom{n-1}{0} \binom{p-1}{0}^{-1} \binom{q-1}{0}^{-1} 1^{-1} 2^{-1} 321 \binom{q-1}{0} \binom{p-1}{0} \binom{n-1}{0}^{-1} 1^{-1} 2^{-1}\\
&&= & \binom{p}{1}^{-1} \binom{q}{1}^{-1}  21 \binom{n-1}{0} 1^{-1} 2^{-1} 321 \binom{n-1}{0}^{-1} 1^{-1} 2^{-1} \binom{q}{1} \binom{p}{1} \\
&&= & \binom{p}{1}^{-1} \binom{q}{1}^{-1}  21 34 1^{-1} 2^{-1} 321 4^{-1} 3^{-1} 1^{-1} 2^{-1} \binom{q}{1} \binom{p}{1} \\
&&= & \binom{p}{1}^{-1} \binom{q}{1}^{-1}  3^{-1} 4^{-1} 1^{-1} 2^{-1} 321 43 \binom{q}{1} \binom{p}{1} \\
&& = & \binom{p}{0}^{-1} \binom{q}{0}^{-1}   1^{-1} 2^{-1} 321  \binom{q}{0} \binom{p}{0} 
\ = \ b_{qp}
\end{array}$$
where in the second last line we used
$$\begin{array}{rcl}
2134 \overline{1} \; \overline{2} 321 \overline{4} \;\overline{3}\; \overline{1} \;\overline{2} = 213 \overline{1} \;\overline{2}43\overline{4} 21 \overline{3} \;\overline{1} \;\overline{2} &=& 213 \overline{1} \;\overline{2} \;\overline{3}43 21 \overline{3} \;\overline{1} \;\overline{2} = 
\overline{3}\; \overline{1}\; \overline{2} 321 4 \overline{1}\; \overline{2} \;\overline{3} 213  \\
& = &\overline{3}\; \overline{1} \;\overline{2} 3 4  \;\overline{3} 213   = \overline{3} \;\overline{1} \;\overline{2}  \;\overline{4} 3 4  213 
= \overline{3}  \; \overline{4} \;\overline{1}\; \overline{2} 3   21 43 ,
\end{array}$$
where the overlining indicates inverse.
\item Finally we consider the generators $b_q^{(i)}$:
 $$\begin{array}{lrcl}  
(i) &\beta b_{n-1}^{(i-2)} \beta^{-1} 
&= &21 \binom{n-1}{0} \binom{n-1}{0}^{-1} b_0^{(i-2)} \binom{n-1}{0} \binom{n-1}{0}^{-1} 1^{-1} 2^{-1}
\ = \ 21  b_0^{(i-2)}  1^{-1} 2^{-1} = b_0^{(i)}\\
\\
(ii) &\beta b_{q-1}^{(i-1)} \beta^{-1} 
&= &21 \binom{n-1}{0} \binom{q-1}{0}^{-1} b_0^{(i-1)} \binom{q-1}{0} \binom{n-1}{0}^{-1} 1^{-1} 2^{-1} \\
&&= &\binom{q}{1}^{-1} 21 \binom{n-1}{0} b_0^{(i-1)} \binom{n-1}{0}^{-1} 1^{-1} 2^{-1} \binom{q}{1} \\
&&= &\binom{q}{0}^{-1} 3 21 3 b_0^{(i-1)} 3^{-1} 1^{-1} 2^{-1} 3^{-1} \binom{q}{0}
\ = \ \binom{q}{0}^{-1} b_0^{(i)} \binom{q}{0} \ = \ b_{q}^{(i)},
\end{array}$$
where in the last line we use $3 21 3 b_0^{(i-1)}  = b_0^{(i)} 3 21 3 $; 
to see this, observe that the defining relations for $B$ imply
$21 (3213) = (3213) 21$ and 
$2^{-1} (3213) = (3213) 1^{-1}$. Firstly, if $i=0$ we have
$$b_0^{(i)} (3213) = 1 (3213) = 2^{-1} 21 (3213) = 2^{-1} (3213) 21 = 
(3213) 1^{-1}21 = (3213) b_0^{(e-1)}.$$
Next, suppose that $i$ is even and larger than 0. Then 
$$\begin{array}{rcl}
\bzi (3213) = \boi (3213) 
&= &\langle 21 \rangle^{(i)} 2^{-1} (3213) \left( \langle 21 \rangle^{i-2} \right)^{-1} \\
&= &\langle 21 \rangle^{(i)} (3213) 1^{-1} \left( \langle 21 \rangle^{i-2} \right)^{-1} \\
&= &(3213) \langle 21 \rangle^{(i)} 1^{-1} \left( \langle 21 \rangle^{i-2} \right)^{-1} \\
&= &(3213) \langle 21 \rangle^{(i-1)} \left( \langle 21 \rangle^{i-2} \right)^{-1}
= (3213) b_0^{(i-1)}.
\end{array}
$$
Finally, suppose that $i$ is odd. Then we have
$$\bzi (3213) 
= 21 \left(b_0^{(i-1)} \right)^{-1} (3213)
= 21 (3213) \left(b_0^{(i-2)} \right)^{-1}
= (3213) 21 \left(b_0^{(i-2)} \right)^{-1}
= (3213) b_0^{(i-1)}.
$$
\end{enumerate}
\qed

\begin{proposition}
\label{bgens} 
The elements $b_{pq}$, $b_{qp}$ and $b_0^{(i)}$ map under $\varphi$ to 
the corresponding generators $a_{pq}$, $a_{qp}$ and $a_0^{(i)}$ 
of $G$ respectively. Thus $\varphi$ is surjective.
\end{proposition}

\proof
Firstly, for each $0 \leq i < e-1$ we have
$b_{i, i+1} = \binom{i+1}{i} \ \phimt \ a_{i,i+1}.$ We 
prove the general case $b_{pq}\  \phimt \ a_{pq}$ by induction on $q-p$:
we use the fact that $\binom{j+1}{i} = \binom{j}{i} \binom{j+1}{j}$,
so that $b_{p, q+1} = \binom{q+1}{q}^{-1} b_{pq} \binom{q+1}{q}$, and hence
$$b_{p, q+1} = \ \phimt \ a_{q, q+1}^{-1} a_{pq} a_{q, q+1} 
=  a_{q, q+1}^{-1} a_{q, q+1}  a_{p, q+1} = a_{p,q+1}.$$
Secondly, we observe that $b_0^{(0)}( = \tau_1)$ and 
$b_0^{(1)}( = \tau_2)$ by definition map to $a_0^{(0)}$ and $a_0^{(1)}$
respectively. Then we continue inductively on $i$:
clearly $b_0^{(i+1)} b_0^{(i)} = 21$ for $i>1$, so 
$$b_0^{(i+1)} = 2 1 \bzii \ \phimt \ a_0^{(1)} a_0^{(0)} ( a_0^{(i)} )^{-1}
= a_0^{(i+1)} a_0^{(i)} ( a_0^{(i)} )^{-1}
= a_0^{(i+1)}.$$
Thus the generators $a_{pq}$ with $p<q$ and $a_0^{(i)}$ 
for all $i$ are in the image of $\varphi$. 
For $p<q$, we have that $a_{qp} = \alpha^{-(n-q)} a_{0 p-q} \alpha^{n-q}
= \varphi(\beta^{-(n-q)} b_{0 p-q} \beta^{n-q})$, so is in $\varphi(B)$.
Finally, writing $X_p$ for  $\{ a_p^{(i)} \mid 0 \leq i <n \}$, we have
$\alpha^{-1} X_p \alpha
= X_{p-1}$ for any $p$, so
$a_p^{(i)}$ will appear in 
$$\varphi \left(\beta^{-(n-p)} \ 
\{ b_0^{(i)} \mid 0 \leq i < n\} \ \beta^{n-p} \right)
= X_p.$$
Thus all the generators of $G$ lie in the image of $\varphi$.
\qed

To show that $B$ and $G$ are isomorphic,
it remains to show that $\varphi$ is one-to-one. This is the case if the 
elements $b_{pq}, b_{qp}$ and $b_{p}^{(i)}$ satisfy the defining relations
for the $a_{pq}, a_{qp}$ and $a_{p}^{(i)}$.

\begin{proposition}
The defining relations for $G$ in terms of $a_{pq}$, etc.,
hold in $B$ in terms of $b_{pq}$, etc.
\end{proposition}

\proof
The defining relations are stated in terms of cyclically ordered
sequences $(p_1,\ldots,p_f)$. Notice that such a sequence
is cyclically ordered precisely when $(p_1-1,\ldots,p_f-1)$ is
cyclically ordered; that is, the action of integer shifting preserves
cyclic order.

Using conjugation by $\beta$, it thus suffices to check the relations for
but one sequence in any orbit under integer shifting: for example, 
suppose that $(p,q,r,s)$ is cyclically 
ordered; then so is $(0,q-p,r-p,s-p)$, and
$b_{pq} b_{rs} = b_{rs} b_{pq}$ if and only if 
$b_{0,q-p} b_{r-p,s-p} = b_{0,q-p} b_{r-p, s-p}$ (conjugation
by $\beta^p$); thus if the relation holds with respect to 
cyclically ordered quadruples of the form $(0,s,t,u)$, 
then it holds for all cyclically ordered quadruples.

Furthermore, the subgroup of $B$ generated
by $b_{i, i+1}$ is isomorphic to the braid group of type $A_n$;
the $b_{pq}$ for $0 \leq p<q<n$ are precisely the 
BKL-generators (see~\cite{bkl}) and
for linearly ordered sequences,
the relations $\R_1$, $\R_2$ and $\R_3$ are precisely the 
BKL-relations these generators; thus they are already known to hold.
Since every cyclically ordered sequence can be integer shifted
to a linearly ordered one, we thus have that $\R_1$, $\R_2$ and 
$\R_3$ hold in $B$.

Now consider any cyclically ordered triple $(p,q,r)$; we want to 
check $\R_4$ holds: that is, for each $i$, 
$b_{pq} b_r^{(i)} = b_r^{(i)} b_{pq}$. Now conjugating by 
$\beta^r$, this holds if and only if
$b_{p-r,q-r} b_{0}^{(j)} = b_0^{(j)} b_{p-r,q-r}$, for some $j$.
Now $0 < p-r < q-r$, and for any $0<s<t$ and any $j$, 
$$\begin{array}{c}
b_{st} b_0^{(j)} = \binom{t}{s+1}^{-1} \binom{t}{s} b_0^{(j)}
= b_0^{(j)} \binom{t}{s+1}^{-1} \binom{t}{s} = b_0^{(j)}b_{st}
\end{array}$$
since $\tau_2$ and $\tau_2'$ commute with $\binom{t}{s}$ whenever
$s \geq 1$.

Consider $0 < q < n$, and any $i$; we will verify $\R_5$ when $p=0$.
Firstly, 
$b_{0q} b_q^{(i)} = \binom{q}{1}^{-1} b_0^{(i)} \binom{q}{0}
=  b_0^{(i)} \binom{q}{1}^{-1} \binom{q}{0}$, since
$\binom{t}{s}$ and $b_0^{(i)}$ commute whenever $s \geq 1$.
This is $b_0^{(i)} b_{0q}$ by definition. On the other hand,
$b_q^{(i)} b_0^{(i)} = \binom{q}{0}^{-1} b_0^{(i)} 3 \binom{q}{1} b_0^{(i)}
= \binom{q}{0}^{-1} b_0^{(i)} 3  b_0^{(i)}\binom{q}{1}$. We will show that
$b_0^{(i)} 3  b_0^{(i)} =  3  b_0^{(i)} 3$ for all $i$; for $i=0,1$
this is precisely the defining relations 
$\tau_2' \tau_3 \tau_2' = \tau_3 \tau_2' \tau_3$ and
$\tau_2 \tau_3 \tau_2 = \tau_3 \tau_2 \tau_3$ respectively.
We continue inductively: let
$2 \leq i+1 < e$. Then
$$\begin{array}{rcl}
 3 {b_0^{(i+1)}}^{-1} 21  3 
\ = \  3 b_0^{(i)} 3            
\ = \  b_0^{(i)} 3  b_0^{(i)}        
\  = \ {b_0^{(i+1)}}^{-1} 21 3  b_0^{(i)} 
&= & {b_0^{(i+1)}}^{-1} 3^{-1} 3 21 3  b_0^{(i)}  \\          
&   = & {b_0^{(i+1)}}^{-1} 3^{-1} b_0^{(i+1)}  3 21 3,
\end{array}$$
where we use $ 3213 b_0^{(i)} = b_0^{(i+1)} 3213$, shown in the 
proof of Lemma~\ref{alpharot} (see page~\pageref{proofalpharot}).
Cancelling $213$ from the right, and rearranging, gives
$ 3 b_0^{(i+1)} 3 = b_0^{(i+1)} 3  b_0^{(i+1)}$ as required.
Thus we have
$b_q^{(i)} b_0^{(i)} = \binom{q}{0}^{-1}  3 b_0^{(i)} 3  \binom{q}{1}
= \binom{q}{1}^{-1}  b_0^{(i)} \binom{q}{0}$, as required.

Now take any $0<p\neq q< n$ and any $i$; a relation of type $\R_5$ 
is of the form $b_{pq} b_q^{(j)} = b_q^{(j)} b_p^{(i)} = b_p^{(i)} b_{pq}$
for a fixed $j$. Conjugating by $\beta^p$, this holds if and only if
$b_{0,q-p} b_{q-p}^{(k)} = b_{q-p}^{(k)} b_0^{(l)} = b_0^{(l)} b_{0,q-p}$
holds, for certain $k$ and $l$. Since the image of the first equation
holds in $G$, its conjugate does as well, and since we know
$b_{0,x} b_{x}^{(k)} = b_{x}^{(k)} b_0^{(k)}$ in $G$, for all $x$ and $k$,
then we have $k=l$. Thus the second equation in terms of the $b$-generators
is of the form shown in the previous paragraph to hold.

Finally, we show that all relations of the form $\R_6$ hold. But for any
$0\leq p < n$ and $0 \leq i < e$, 
$b_p^{(i)} b_p^{(i-1)} = \binom{p}{0}^{-1} b_0^{(i)} b_0^{(i-1)} \binom{p}{0}$.
We observed earlier that $b_0^{(i)} b_0^{(i-1)} = 21$ for all $0<i<e$; 
and if $i = 0$, we have
$
b_0^{(0)} b_0^{(e-1)} 
= 1 \langle 21 \rangle^{e-1} \left(\langle 21 \rangle^{e-2} \right)^{-1} 
=  \langle 12 \rangle^{e} \left(\langle 21 \rangle^{e-2} \right)^{-1}
=  \langle 21 \rangle^{e} \left(\langle 21 \rangle^{e-2} \right)^{-1} = 21.
$
Thus $b_p^{(i)} b_p^{(i-1)} = \binom{p}{0}^{-1} 21 \binom{p}{0}$, which
is independent of $i$, so $\R_6$ holds. 
\qed

Thus we can complete the main result.

\begin{theorem}
\label{main}
The presentation of Section~\ref{statement} gives rise to a Garside structure for
the braid group $B$ of the (complex) reflection group $G(e,e,r)$, 
where the Garside element is the
least common multiple of the generators in the Garside monoid; and furthermore,
for which the generators map to reflections under the  natural map 
$\nu:B \rightarrow G(e,e,r)$.\end{theorem}

\proof
The only part of the assertion remaining unproved is the last phrase. We know 
from~\cite{bmr} that the images of the BMR-generators under $\nu$ are reflections
in $G(e,e,r)$. From the definitions given immediately before the statement of 
Lemma~\ref{alpharot}, we see that each of the new generators is a conjugate
of a BMR-generator. Thus its image under $\nu$ is a conjugate of a reflection,
and hence is a reflection.
\qed

\section{No finite nice presentation on BMR-generators exists}

\noindent
We will call a presentation of a braid group ``nice'' if it is of the form
given in the above result.

In this section, 
we will show that there is no finite presentation of the braid group
of type $G(e,e,r)$ on the BMR-generators, 
where the relations are all positive words over these generators, such
that the corresponding monoid embeds in the braid group.
Recall that we have an infinite collection of equations in the braid group
of the form $21^n3213 = 32132^n1$. In the monoid, however,
none of these implies another, and for all $n>1$, none of these is implied 
by the defining relations.
Thus it suffices to show that there is no finite collection of relations
on the BMR-generators of the required form which imply all of these equations.

In order to derive $21^n3213 = 32132^n1$ we need to be able to rewrite
$21^n3213$; thus we need to rewrite some word of the form $21^k$ or
of the form $1^k 3213$ in terms of the generators $\{1,2,3\}$, where
``rewrite $w_1$ in terms of the generators $\{1,2,3\}$'' means: find an equation 
$w_1 = w_2$ in the braid group $B$ where $w_2$ is a (positive) word over 
$\{1,2,3\}$. Let $M$ denote the monoid generated by the  $b_{pq}, b_{qp}$ and
$b_p^{(i)}$; since $M$ embeds in $B$ this is exactly the same as 
finding an equation $w_1 = w_2$ which holds in $M$ where $w_2$ 
is a word over $\{1,2,3\}$. Since $M$ is a Garside monoid,
we have solutions to the word problem and can calculate lcm's,
which will do the work for us here. Recall that
$1 = a_0^{(0)}, 2 = a_0^{(1)}$ and $3 = a_{01}$.

\begin{proposition}
\label{nofinite}
There is no finite presentation for the braid group of type $G(e,e,r)$
with generating set that given in~\cite{bmr}
and with all relations positive words, for which the corresponding monoid
embeds in the braid group.
\end{proposition}

\proof
There is a surjection $B(e,e,n+1) \rightarrow B_n$, the braid group of
type $A$ on $n$ strings, given by $1,2 \mapsto \sigma_1$ and
 $i \mapsto \sigma_{i-1}$ for $i \geq 3$. Thus 
$21^k \stackrel{\varphi}{\mapsto} \sigma_1^{k+1}$; so $2 1^k$ 
can only be rewritten in terms of $\{2,1\}$. We may suppose that 
$21^k = 1 w$ for some word $w$ over $\{2,1\}$; we have that
$21^k = 21  1^{k-1} = 1 b_0^{(e-1)} \ 1^{k-1}$; so by left cancellation,
$w = b_0^{(e-1)} \ 1^{k-1}$. However 
$b_0^{(e-1)} \ 1^{k-1} \equiv b_0^{e-1} \ (b_0^{(0)})^{k-1}$
is in a singleton equivalence class in $M$, so can never be 
rewritten in terms of $\{1,2\}$. Thus it remains to 
show that $1^k 3213$ cannot be rewritten in terms
of $\{1,2,3\}$. 

Let $U$ denote the word $1^k 3213$; it suffices to show that
for any $k \geq 0$, we cannot rewrite $U$ in terms of $\{1,2,3\}$.
Since $M$ is Garside, then it has the `reduction property', and
so we can quickly determine divisibility in $M$ by a generator using 
the method of $a$-chains (see~\cite{corran}). 

Observe that $U = 3213 u$ where $u = (a_0^{(e-2)})^k$; clearly
$u$ is not divisible by $1$, $2$, or $3$, so is certainly not 
rewriteable over these letters. Also, neither $1$ nor $2$ divides
$3u$ so $U$ cannot be rewritten over $\{1,2,3\}$ to begin with $321$.
Since $13 = 3 a_1^{(0)}$, we have $U = 323 v$ where $v = a_{1}^{(0)} u$,
which is not divisible by $1, 2$ or $3$; furthermore, $3v$ is not divisible
by $2$, so we see that $U$ cannot be rewritten $\{1,2,3\}$ to begin with $32$.
Next, we write $U = 313 w $ where $w = a_1^{(e-1)} u$; this
is not divisible by 1,2 or 3; and $3w$ is not divisible by 1 or 2. Thus
$U$ cannot be rewritten over $\{1,2,3\}$ to start with  $31$.
Finally, $U = 3^2 x$ where $x = a_1^{(1)} a_1^{(0)} u$;
$x$ is not divisible by 1,2 or 3, so we have that $U$ cannot be rewritten
to begin with $33$ either. This exhausts all the possibilites for rewriting
$U$ over $\{1,2,3\}$ starting with $3$.

On the other hand, we can write $U = 213y$ where $y= a_{10} u$, which is not
divisible by $1,2$ or $3$; so not rewriteable over these letters. Moreover,
$3y$ is not divisible by $1$ or $2$; and $13y$ is not divisible by $2$.
Thus we have that $U$ cannot be rewritten to begin with $21$ or $22$. 
So now we write $U = 232 v$ where $v = a_{1}^{(0)} u$, 
which is not divisible by $1, 2$ or $3$; furthermore, $2v$ is not divisible
by $1$ or $3$. In this way we have exhausted all the possibilites for rewriting
$U$ over $\{1,2,3\}$ starting with $2$.

Finally, we want to show that $U$ cannot be rewritten to begin with $1$. By 
left cancellation, we may assume that $k=0$; then we have that
$3213 = 131a_1^{-1}$. But $31a_1^{-1}$ (being in a singleton equivalence class
in $M$) cannot be rewritten over $\{1,2,3\}$. The result follows.
\qed


\begin{tabular}{l}
{\sc David Bessis} \\
LIFR-MI2P, CNRS / Independent University of Moscow\\
Bolshoi Vlasevsky Pereulok, Dom 11, Moscow 121002, Russia\\
Email:  bessis@mccme.ru\\
\\
{\sc Ruth Corran}\\
Institut Henri Poincar\'e\\
11 rue Pierre et Marie Curie,
75005 Paris, France\\
Email: corran@math.jussieu.fr
\end{tabular}

\end{document}